\documentclass[11pt]{article}
\usepackage{amsfonts}
\usepackage{bbm}
\usepackage{amscd}
\usepackage{xypic}
\usepackage{amsmath}
\usepackage{amssymb}
\usepackage{amsthm}
\usepackage{bbm}
\usepackage{CJK}
\usepackage{fancyhdr}
\usepackage{hyperref}
\usepackage{indentfirst}
\usepackage{latexsym}
\usepackage{mathrsfs}
\allowdisplaybreaks 
\usepackage{color}

\usepackage[top=1in,bottom=1in,left=1.25in,right=1.25in]{geometry}
\def\<{\langle}
\def\>{\rangle}
\def\a{\alpha}
\def\b{\beta}
\def\ci{\circ}
\def\c{\cdot}

\def\D{\Delta}

\def\g{\gamma}

\def\lr{\longrightarrow}

\def\o{\otimes}

\def\v{\varepsilon}

\def\<{\langle}
\def\>{\rangle}

\date{}
\begin{document}
	\renewcommand{\baselinestretch}{1.2}
	\renewcommand{\arraystretch}{1.0}
	\title{\bf Radford $[(m,k),m]$-biproduct Theorem for Generalized Hom-crossed Products}
	\date{}
\author {{\bf Botong Gai}$^{(1)}$\,and \, {\bf Shuanhong Wang$^{(2)}$\footnote {Corresponding author:  shuanhwang@seu.edu.cn}}\\
{\small (1): School of Mathematics, Southeast University,  Nanjing, Jiangsu}\\
{\small 210096, P. R. of China. E-mail: 230228425@seu.edu.cn}\\
{\small (2):  Shing-Tung Yau Center, School of Mathematics, Southeast University}\\	
{\small Nanjing, Jiangsu 210096, P. R. of China}}
 \maketitle
	\begin{center}
		\begin{minipage}{14.cm}
			{\bf Abstract}  In this paper, we mainly provide a new approache to construct Hom-Hopf algebras. For this, we introduce and study
 the notion of a left $(m,k)$-Hom-crossed product structure as a generalization of  $k$-Hom-smash product structure. Then one  combines this
  $(m,k)$-Hom-crossed product structure and a left
  $m$-Hom-smash coproduct structure to build Radford $[(m,k),m]$-biproduct theorem.  Finally, we study Hom admissible mappping system
   to characterize  this Radford $[(m,k),m]$-biproduct structure.
	\\
		
			\\

		\end{minipage}
	\end{center}

\section*{Introduction}
	\def\theequation{0. \arabic{equation}}
	\setcounter{equation} {0} \hskip\parindent
	
In the classical Hopf algebraic theory,  the one of the celebrated
 results  is  Radford' biproduct \cite{R85} which provided in particular
 an important approach to solve the classification of finite-dimensional pointed
  Hopf algebras (see \cite{AS98} and \cite{AS10}).  This biproduct says that if $A$ is a braided Hopf algebra in the braided monoidal category of Yetter-Drinfeld modules ${}^H_H{\cal YD}$ over a Hopf algbera $H$, then a left smash product algebra structure and a left smash coproduct coalgebra structure
   afford a Hopf algebra structure  on  $A\o H$, see \cite{M94}.
   Radford's biproduct was generalized to many cases: replacing  the smash product algebra structure by  a left Hopf crossed product (see, \cite{WJZ})
    and replacing  the smash coproduct coalgebra structure by  a left Hopf crossed coproduct  \cite{JWZ}  in the setting of Hopf algebras; replacing Hopf algebras by quasi-Hopf algebras \cite{BN}, by multiplier Hopf algebras \cite{D} and by monoidal Hom-Hopf algebras (see, \cite{LS14},\cite{MLC}).
\\

As we know that  the notion of a left Hopf crossed product was  introduced in \cite{BCM} and the dual Hopf crossed coproduct was
  introduced in \cite {DRZ} (see, \cite{W1,W2}).  These notions have been
  studied in the setting of  Hom-Hopf algebras (see, \cite{LW14,LW19}).
\\

The main object of this paper is to provide  new methods to construct
 monoidal Hom-Hopf algebras as a generalization of the one in  \cite{WJZ} by  introducing and studying the notions of a $(m,k)$-Hom-crossed product.
 \\

The organization of the paper is the following. In Section 1, some basic notations
 about monoidal Hom-Hopf algebras, Hom-(co)module algebras,  Hom-(co)module coalgebras and other materials that we will need are recalled.
 In Section 2 and Section 3, we will introduce the notion of a left $(m,k)$-Hom-crossed product for a monoidal Hom-Hopf algebra
 and obtain  Radford $[(m,k),m]$-biproduct with $ m, k\in \mathbb{Z}$ (see Theorem 2.3, Theorem 3.3 and Theorem 3.5). In Section 4,
  in order to characterize the Radford $[(m,k),m]$-biproduct structure, we study weak Hom admissible mappping system (see, Theorem 4.11).
   \\

Throughout, let $k$ be a fixed field and everything is over $k$ unless otherwise
  specified.  We refer the readers to the book of Sweedler \cite{S}
	for the relevant concepts on the general theory of Hopf
	algebras.  Let $(C, \Delta )$ be a coalgebra, we use the Sweedler-Heyneman's
   notation for $\Delta $ as follows: $\Delta (c)=\sum c_1\otimes c_2$, for all $c\in C$.

\section*{1. Preliminaries}
	\def\theequation{1. \arabic{equation}}
	\setcounter{equation} {0} \hskip\parindent

In this section we will recall the notions of  a monoidal Hom-category,
  a monoidal Hom-Hopf algebra,
	 a (co)action of monoidal Hom-Hopf algebra and a Hom-smash (co)product.

\subsection* {1.1. A monoidal Hom-category $\tilde{\mathcal{H}}(\mathcal {M}_{k})$}

	Let $\mathcal{M}_{k}=(\mathcal{M}_{k},\o,k,a,l,r )$
	denote the usual monoidal category of $k$-vector spaces and linear maps between them.
	Recall from [6]
	that there is the {\it monoidal Hom-category} $\widetilde{\mathcal{H}}(\mathcal{M}_{k})=
	(\mathcal{H}(\mathcal{M}_{k}),\,\o,\,(k,\,id),
	\,\widetilde{a},\,\widetilde{l},\,\widetilde{r })$, a new monoidal category,
	associated with $\mathcal {M}_{k}$ as follows:
	
	$\bullet$
	$ \mathcal{H}(\mathcal{M}_{k})$ are couples
	$(M,\mu)$, where $M \in \mathcal {M}_{k}$ and $\mu \in Aut_k(M)$, the set of
	all $k$-linear automorphisms of $M$;
	
	$\bullet$   $f:(M,\mu)\rightarrow (N,\nu)$ in $ \mathcal{H}(\mathcal{M}_{k})$
	is a $k$-linear map $f: M\rightarrow N$ in $\mathcal{M}_{k}$
	satisfying   $ \nu \circ f = f\ci \mu$  for any two objects
	$(M,\mu),(N,\nu)\in \mathcal{H}(\mathcal{M}_{k})$;
	
	$\bullet$   The tensor product is given by

\begin{equation*}
(M,\mu)\otimes (N,\nu)=(M\otimes N,\mu\otimes\nu )
\end{equation*}
	for any $(M,\mu),(N,\nu)\in \mathcal{H}(\mathcal{M}_{k})$;
	
	$\bullet$ The tensor unit is given by $(k, id)$;
	
	$\bullet$   The associativity constraint $\widetilde{a}$
	is given by the formula

\begin{equation*}
\widetilde{a}_{M,N,L}=a_{M,N,L}\circ((\mu\o id)\o
	\varsigma^{-1})=(\mu\o(id\o\varsigma^{-1}))\circ a_{M,N,L},
\end{equation*}
	for any objects $(M,\mu),(N,\nu),(L,\varsigma)\in \mathcal{H}(\mathcal{M}_{k})$;
	
	$ \bullet $  The left and right unit constraint
	$\widetilde{l}$ and $\widetilde{r }$ are given by

\begin{equation*}
\widetilde{l}_M=\mu\circ l_M=l_M\circ(id\otimes\mu), \quad
	\widetilde{r}_M =\mu\circ r_M=r_M\circ(\mu\otimes id)
\end{equation*}
	for all $(M,\mu) \in \mathcal{H}(\mathcal{M}_{k})$.
	\\

\subsection* {1.2. Monoidal Hom-associative algebras and
  monoidal Hom-coassociative coalgebras}

A {\it unital monoidal Hom-associative algebra} is a vector space $A$
	together with an element $1_A\in A$ and linear maps
\begin{equation*}
m: A\otimes A\lr A (a\otimes b\mapsto ab) \quad \mbox{and} \quad \alpha\in Aut_k(A)
\end{equation*}
	such that
	\begin{equation}
		\alpha(a)(bc)=(ab)\alpha(c),
	\end{equation}
	\begin{equation}
		\alpha(ab)=\alpha(a)\alpha(b),
	\end{equation}
	\begin{equation}
		a1_A=1_Aa=\alpha(a),
	\end{equation}
	\begin{equation}
		\alpha(1_A)=1_A
	\end{equation}
	for all $a,b,c\in A.$
\\
	
	{\bf Remark.}   (1) In the language of algebras, $m$ is called the
	Hom-multiplication, $\alpha$ is the twisting
	automorphism and $1_A$ is the unit. Note that Eq.(1.1) can be rewritten as
	$a(b\alpha^{-1}(c)) = (\alpha^{-1}(a)b)c$.
	The monoidal Hom-algebra $A$ with a {\sl structure map} $\alpha$ will be denoted by $(A,\alpha)$.

(2) A monoidal Hom-associative algebra is not the same as a Hom-associative algebra
 in which $\a $ is not necessary bijective,
 (see, \cite{MP} and \cite{MS}).
	
	(3)  Let $(A,\alpha)$ and  $(A',\alpha')$ be two monoidal Hom-algebras.
	A monoidal Hom-algebra map
	$f:(A,\alpha)\rightarrow (A',\alpha')$ is a linear map such that
	$f\circ \alpha=\alpha'\circ f, f(ab)=f(a)f(b)$ and
	$f(1_A)=1_{A'}.$
	\\	
	
	{\it A counital monoidal Hom-coassociative coalgebra} is
	a vector space $C$ together with linear maps
	$\D:C\lr C\otimes C (\Delta(c)=c_1\otimes c_2)$,
	$\varepsilon:C \lr k$ and $\g \in Aut_k(C)$ so that
	\begin{equation}
		\gamma^{-1}(c_1)\otimes\Delta(c_2)=\Delta(c_1)\otimes\gamma^{-1}(c_2),
	\end{equation}
	\begin{equation}
		\Delta(\gamma(c))=\gamma(c_1)\otimes\gamma(c_2),
	\end{equation}
	\begin{equation}
		c_1\varepsilon(c_2)=\gamma^{-1}(c)=\varepsilon(c_1)c_2,
	\end{equation}
	\begin{equation}
		\varepsilon(\gamma(c))=\varepsilon(c)
	\end{equation}
	for all $c\in C.$
	\\
	
	{\bf Remark.}   (1) Note that Eq.(1.5) is equivalent to
	$c_1\otimes c_{21}\otimes \gamma(c_{22})=\gamma(c_{11})\otimes c_{12}\otimes c_2.$
	Similar to monoidal Hom-algebras, monoidal Hom-coalgebras
	will be short for counital monoidal Hom-coassociative coalgebras
	without any confusion.  The monoidal Hom-coalgebra $C$ with a {\sl structure map} $\g $ will be denoted by $(C,\g )$.

(2) A monoidal Hom-coassociative coalgebra is not the same as a Hom-coassociative coalgebra in which
 Eq.(1.5) and Eq.(1.7) are replaced by $\gamma (c_1)\otimes\Delta(c_2)=\Delta(c_1)\otimes\gamma (c_2)$
  and  $c_1\varepsilon(c_2)=\gamma (c)=\varepsilon(c_1)c_2$ for any $C\in C$, respectively, and
   $\g $ is not necessary bijective, (see, \cite{MP} and \cite{MS}).
	
	(3)  Let $(C,\gamma)$ and $(C',\gamma')$ be two monoidal Hom-coalgebras.
	A monoidal Hom-coalgebra map $f:(C,\gamma)\rightarrow(C',\gamma')$
	is a linear map such that $f\circ \gamma=\gamma'\circ f, ~\Delta\circ f=(f\otimes f)\circ\Delta$
	and $\varepsilon'\circ f=\varepsilon.$

\subsection* {1.3. Monoidal Hom-Hopf algebras}
	
	{\it A monoidal Hom-bialgebra} $H=(H,\alpha,m,1_H,\D,\varepsilon)$
	is a bialgebra in the monoidal category
	$ \tilde{\mathcal{H}}(\mathcal {M}_{k}).$
	This means that $(H,\alpha,m,1_H)$ is a monoidal Hom-algebra and
	$(H,\alpha,\D,\varepsilon)$ is
	a monoidal Hom-coalgebra such that $\D$ and $\varepsilon$
	are morphisms of algebras.

	A monoidal Hom-bialgebra $(H, \alpha)$ is called {\it a monoidal Hom-Hopf algebra}
	if there exists a morphism (called {\sl antipode})
	$S: H\rightarrow H$ in $ \tilde{\mathcal{H}}(\mathcal {M}_{k})$
	(i.e. $S\ci \alpha=\alpha\ci S$),
	which is the convolution inverse of the identity morphism $id_H$
	(i.e. $ S*id=1_H\ci \varepsilon=id*S $). Explicitly,  for all $h\in H$,

\begin{equation*}
S(h_1)h_2=\varepsilon(h)1_H=h_1S(h_2).
\end{equation*}
\\
	
	{\bf Remark.}   (1) Note that a monoidal Hom-Hopf algebra is
	by definition a Hopf algebra in the category $ \tilde{\mathcal{H}}(\mathcal {M}_{k})$
 (see, \cite{CG}).
	
	(2)  Furthermore, the antipode of monoidal Hom-Hopf algebras has
	almost all the properties of antipode of Hopf algebras such as

\begin{equation*}
S(hg)=S(g)S(h), S(1_H)=1_H,
	\D(S(h))=S(h_2)\o S(h_1),\, \, \mbox{and}\, \, \varepsilon\ci S=\varepsilon.
\end{equation*}
	That is, $S$ is a monoidal Hom-anti-(co)algebra homomorphism.
	Since $\alpha$ is bijective and commutes with $S$,
	we can also have that the inverse $\alpha^{-1}$ commutes with $S$,
	that is, $S\ci \alpha^{-1}= \alpha^{-1}\ci S.$

(3) A monoidal Hom-Hopf algebra $(H, \a )$ is not the same as a Hom-Hopf algebra $(H, \a )$
 in which $\a $ is not necessary bijective,
 (see, \cite{MP} and \cite{MS}).

\subsection* {1.4. $m$-Hom-smash products and $m$-Hom-smash coproducts}

Let $(A,\alpha)$ be a monoidal Hom-algebra.
	{\it A left $(A,\alpha)$-Hom-module} consists of
	an object $(M,\mu)$ in $\tilde{\mathcal{H}}(\mathcal {M}_{k})$
	together with a morphism
	$\psi:A\o M\lr M,~\psi(a\o m)=a\cdot m$ such that

\begin{equation*}
\alpha(a)\c(b\c m)=(ab)\c\mu(m),\,\,
	\,\,\mu(a\c m)=\alpha(a)\c\mu(m),\,\, \mbox{and}
	\,\,1_A\c m=\mu(m),
\end{equation*}
	for all $a,b\in A$ and $m \in M$.
\\
	
	{\bf Remark.} (1) Monoidal Hom-algebra $(A,\alpha)$ can be
	considered as a Hom-module on itself by the Hom-multiplication.

(2) Let $(M,\mu)$ and $(N,\nu)$ be two left $(A,\alpha)$-Hom-modules.
	A morphism $f: M\lr N$ is called a left
	$(A,\alpha)$-linear if
	$f(a\c m)=a\c f(m),~f\ci \mu= \nu\ci f$.
	The category of left $(A,\alpha)$-Hom-modules is denoted by
	$\tilde{\mathcal{H}}(_{A}\mathcal {M})$.
	\\
	
	Similarly, let $(C,\gamma)$ be a monoidal Hom-coalgebra.
	{\it A right $(C,\gamma)$-Hom-comodule} is an object
	$(M,\mu)$ in $\tilde{\mathcal{H}}(\mathcal {M}_{k})$
	together with a $k$-linear map
	$\rho_M: M\lr M\o C,~\rho_M(m)=m_{(0)}\o m_{(1)}$ such that
	\begin{equation}
		\mu^{-1}(m_{(0)})\o \D_C(m_{(1)})
		=(m_{(0)(0)}\o m_{(0)(1)})\o \gamma^{-1}(m_{(1)}),
	\end{equation}
	\begin{equation}
		\rho_M(\mu(m))=\mu(m_{(0)})\o\gamma(m_{(1)}),
		\ \ \
		m_{(0)}\varepsilon(m_{(1)})=\mu^{-1}(m),
	\end{equation}
	for all $m\in M.$
\\	
	
	{\bf Remark.} (1) $(C,\gamma)$ is a Hom-comodule on itself via the Hom-comultiplication.
	
(2) Let $(M,\mu)$ and $(N,\nu)$ be two right $(C,\gamma)$-Hom-comodules.
	A morphism $g: M\lr N$ is called right $(C,\gamma)$-colinear
	if $g\ci \mu=\nu\ci g$ and
	$g(m_{(0)})\o m_{(1)}=g(m)_{(0)}\o g(m)_{(1)}.$
	The category of right
	$(C,\gamma)$-Hom-comodules is denoted by
	$\tilde{\cal{H}}(\cal {M}^C)$ .
	\\
	
	Let $(H, \alpha)$ be a monoidal Hom-bialgebra. A monoidal Hom-algebra $(B,\beta)$ is called a left {\it
		weak $(H,\alpha)$-Hom-module algebra}, if $(B,\beta)$ has a map: $\cdot : H\o B\lr B$,
  obeying the following axioms:
	\begin{equation}
		h\c(ab)=(h_1\c a)(h_2\c b),\,\,\,\,\,\,\,\,h\c1_B=\varepsilon(h)1_B,
	\end{equation}
	for all $a,b\in B,h\in H.$ Furthermore, $(B,\beta)$ is called {\it a left
		$(H,\alpha)$-Hom-module algebra} if $(B,\beta)$ is a left
	$(H,\alpha)$-Hom-module with the action "$\cdot$".

	Let $(H, \alpha)$ be a monoidal Hom-bialgebra. A monoidal Hom-coalgebra $(B,\beta)$ is called {\it a left
		$(H,\alpha)$-Hom-comodule coalgebra}, if $(B,\beta)$ is a left
	$(H,\alpha)$-Hom-comodule with coaction $\rho$ obeying the following axioms:
	\begin{equation}
b_{(-1)}\otimes b_{(0)1}\otimes b_{(0)2}=b_{1(-1)}b_{2(-1)}\otimes b_{1(0)}\otimes b_{2(0)},~\varepsilon(b_{(0)})b_{(-1)}=\varepsilon(b)1_{H}.
	\end{equation}
    for all $ b\in B $.

\section*{2. $(m,k)$-Hom-crossed products}
	\def\theequation{2. \arabic{equation}}
	\setcounter{equation} {0} \hskip\parindent

Let $m, k\in \mathbb{Z}$.  In order to obtain  Radford $[(m,k),m]$-biproduct theorem, in this section, we introduce
 and study the notion of a left $(m,k)$-Hom-crossed product for a monoidal Hom-Hopf algebra.

	{\bf Definition 2.1.} Let $ (H,\alpha) $ be a monoidal Hom-Hopf algebra. Let $ (A,\beta) $ be a left weak $(H,\alpha) $-Hom-module algebra
  and $\sigma: H\otimes H\lr A $ is a linear map. Then
   $(A\sharp _{\sigma}H,\beta\otimes\alpha) $ is called a  {\sl left $(m,k) $-Hom-crossed product} with $m, k\in \mathbb{Z} $ if
	
	(1) $ A\sharp_{\sigma} H=A\otimes H $ as a vector space, and
	
	(2) it has a Hom-multiplication:
		\begin{align*}
			(a\sharp_{\sigma} h)(b\sharp_{\sigma} g)=a[(\alpha^{m}(h_{11})\cdot\beta^{-2}(b))\sigma(\alpha^{k+1}(h_{12}),\alpha^{k}(g_{1}))]\sharp_{\sigma}\alpha(h_{2}g_{2})
		\end{align*}
 for all $a\sharp_{\sigma} h$, $b\sharp_{\sigma} g\in A\sharp_{\sigma} H$.	
\\

	{\bf Remark 2.2.}  (1)	When $ \alpha=id_{H}$ and $ \beta=id_{A} $, we will obtain the usual Hopf crossed product
   (see,\cite{BCM}, \cite{WJZ});

 (2)  When $ m=0$ and  $ k=-1 $, we obtain the monoidal Hom-crossed product (see, \cite{LW14}).
	
 (3) Consider the case when $ \sigma $ is trivial, that is, $ \sigma(h,g)=\varepsilon_{H}(h)\varepsilon_{H}(g)1_{A} $ for all $ h, g\in H $. Then the $ (m,k) $-Hom-crossed product is reduced to $ m $-Hom-smash product \cite{MLC}. In fact, we have
 \begin{align*}
(a\sharp_{\sigma} h)(b\sharp_{\sigma} g)=&a[(\alpha^{m}(h_{11})\cdot\beta^{-2}(b))\sigma(\alpha^{k+1}(h_{12}),\alpha^{k}(g_{1}))]\sharp_{\sigma}\alpha(h_{2}g_{2})
\\=&a[(\alpha^{m}(h_{11})\cdot\beta^{-2}(b))\varepsilon(h_{12})\varepsilon(g_{1})1_{A}]\sharp_{\sigma}\alpha(h_{2}g_{2})
\\=&a(\alpha^{m}(h_{1})\cdot\beta^{-1}(b))\sharp_{\sigma}\alpha(h_{2})g.
\end{align*}
	\\
	
	{\bf Theorem 2.3.} Let $ (A\sharp_{\sigma} H,\beta\otimes\alpha) $ be a $ (m,k) $-Hom-crossed product. Then $ (A\sharp_{\sigma} H,\beta\otimes\alpha) $
    is a monoidal Hom-algebra with unit $ 1_{A}\sharp_{\sigma}1_{H} $ if and only if the following equations hold:
	
	(1) $ \sigma(h,1_{H})=\sigma(1_{H},h)=\varepsilon(h)1_{A} $
	, $ \sigma\circ(\alpha\otimes\alpha)=\beta\circ\sigma $;
	
	(2) $ [\alpha^{m}(h_{1}l_{1})\cdot a]\sigma(\alpha^{k+2}(h_{2}),\alpha^{k+2}(l_{2}))=\sigma(\alpha^{k+2}(h_{1}),\alpha^{k+2}(l_{1}))[\alpha^{m}(h_{2}l_{2})\cdot a] $;
	
	(3) $ [\alpha^{m+1}(h_{1})\cdot\sigma(\alpha^{k+1}(l_{1}),\alpha^{k+1}(g_{1}))]\sigma(\alpha^{k+2}(h_{2}),\alpha^{k+1}(l_{2}g_{2}))$
	 $=\sigma(\alpha^{k+2}(h_{1}),\alpha^{k+2}(l_{1}))\\\sigma(\alpha^{k+1}(h_{2}l_{2}),
 \alpha^{k+1}(g)) $,
for any $h, l, g\in H$ and $a\in A$.

    {\bf Proof. }  We firstly prove the direction "$ \Rightarrow $". If $ (A\sharp_{\sigma} H,\beta\otimes\alpha) $ is a monoidal Hom-Hopf algebra, then
    \begin{align*}
    	 (\beta\otimes\alpha)[(1_{H}\sharp_{\sigma}h)(1_{H}\sharp_{\sigma}g)]=(1_{H}\sharp_{\sigma}\alpha(h))(1_{H}\sharp_{\sigma}\alpha(g)).
    \end{align*}
    Applying $ id_{A}\otimes\varepsilon_{H} $ to both sides, we have
    \begin{align*}
    	 \beta^{3}\circ\sigma(\alpha^{k-1}(h),\alpha^{k-1}(g))=\beta^{2}\circ\sigma(\alpha^{k}(h),\alpha^{k}(g)).
    \end{align*}
    Since $ \alpha$, $\beta $ are bijection, we obtain $ \sigma\circ(\alpha\otimes\alpha)=\beta\circ\sigma $.

    Meanwhile,
    \begin{align*}
    	1_{A}\sharp _{\sigma}\alpha(h)=(1_{A}\sharp _{\sigma} h)( 1_{A}\sharp _{\sigma} 1_{H})=\sigma(\alpha^{k+2}(h_{1}),1_{H})\sharp _{\sigma}\alpha^{2}(h_{2}).
    \end{align*}
Applying $ id_{A}\otimes\varepsilon_{H}  $ to both sides, we have $ \sigma(\alpha^{k+1}(h),1_{H})=\varepsilon_{H}(h)1_{A} $ and $ \sigma(1_{H},\alpha^{k+1}(h))=\varepsilon_{H}(h)1_{A} $. Since $ \alpha $ is bijective and $ \varepsilon_{H}\circ\alpha=\varepsilon_{H} $, then $ \sigma(h,1_{H})=\varepsilon_{H}(h)1_{A}=\sigma(1_{H},h) $.

Besides,
\begin{equation*}
(1_{A}\sharp _{\sigma}\alpha(h))[(1_{A}\sharp _{\sigma} l)(a\sharp _{\sigma} 1_{H})]=[(1_{A}\sharp _{\sigma} h)(1_{A}\sharp _{\sigma} l)](\beta(a)\sharp _{\sigma}{1_{H}}).
\end{equation*}

Applying $ id_{A}\otimes \varepsilon_{H}  $ to both sides, we will obtain Eq.(2).

\begin{equation*}
(1_{A}\sharp _{\sigma}\alpha(h))[(1_{A}\sharp _{\sigma} l)(1_{A}\sharp _{\sigma} g)]=[(1_{A}\sharp _{\sigma} h)(1_{A}\sharp _{\sigma} l)](1_{A}\sharp _{\sigma}\alpha(g))  .
\end{equation*}
Applying $ id_{A}\otimes \varepsilon_{H}  $ to both sides, we will obtain Eq.(3).

Conversely, if (1)-(3) hold. Then for all $ a$, $b$, $c\in A$, $h$, $g$, $l\in H  $, it is easy to verify the following equation:
	\begin{center}
		$ (\beta\otimes\alpha)(1_{A}\sharp _{\sigma}{1_{H}})=(1_{A}\sharp _{\sigma}{1_{H}}), $
	
	$(1_{A}\sharp _{\sigma}{1_{H}})(a\sharp _{\sigma}h)=\beta(a)\sharp _{\sigma}\alpha(h)=(a\sharp _{\sigma} h)(1_{A}\sharp _{\sigma}{1_{H}})$,
	
   $(\beta\otimes\alpha)((a\sharp _{\sigma} h)(b\sharp _{\sigma} g))=(\beta(a)\sharp _{\sigma}\alpha(h))(\beta(b)\sharp _{\sigma}\alpha(g)) $ .
	\end{center}

	We also have
	
	$ [(a\sharp _{\sigma} h)(b\sharp _{\sigma} g)][\beta(c)\sharp _{\sigma}\alpha(l)]\\=[a[(\alpha^{m}(h_{11})\cdot\beta^{-2}(b))\sigma(\alpha^{k+1}(h_{12}),\alpha^{k}(g_{1}))]\sharp _{\sigma}\alpha(h_{2}g_{2})][\beta(c)\sharp _{\sigma}\alpha(l)]
	 \\=[a(\alpha^{m+1}(h_{11})\cdot\beta^{-1}(b))][[\sigma(\alpha^{k+1}(h_{12}),\alpha^{k}(g_{1}))(\alpha^{m}(h_{211}g_{211})\cdot\beta^{-2}(c))]\\~~~~~\sigma(\alpha^{k+2}(h_{212}g_{212})),\alpha^{k+1}(l_{1})]\sharp _{\sigma}\alpha^{2}(h_{22}g_{22})\alpha^{2}(l_{2})
	 \\=[a(\alpha^{m+1}(h_{11})\cdot\beta^{-1}(b))][[\sigma(\alpha^{k+3}(h_{1211}),\alpha^{k+2}(g_{111}))(\alpha^{m+1}(h_{1211})\alpha^{m}(g_{112})\cdot\beta^{-2}(c))]\\~~~~~\sigma(\alpha^{k+2}(h_{112})\alpha^{k+1}(g_{12})),\alpha^{k+1}(l_{1})]\sharp _{\sigma}\alpha(h_{2}g_{2})\alpha^{2}(l_{2})
	 \\\overset{(2)}{=}[a(\alpha^{m+1}(h_{11})\cdot\beta^{-1}(b))][[[\alpha^{m}(\alpha(h_{1211})g_{111})\cdot\beta^{-2}(c)]\sigma(\alpha^{k+3}(h_{1212}),\alpha^{k+2}(g_{112}))]\\~~~~~\sigma(\alpha^{k+2}(h_{122})\alpha^{k+1}(g_{12})),\alpha^{k+1}(l_{1})]\sharp _{\sigma}\alpha(h_{2}g_{2})\alpha^{2}(l_{2})
	 \\=[a(\alpha^{m+1}(h_{11})\cdot\beta^{-1}(b))][[[\alpha^{m-1}(h_{12})\alpha^{m}(g_{111})\cdot\beta^{-2}(c)]\sigma(\alpha^{k+2}(h_{211}),\alpha^{k+2}(g_{112}))]\\~~~~~\sigma(\alpha^{k+2}(h_{212})\alpha^{k+1}(g_{12})),\alpha^{k+1}(l_{1})]\sharp _{\sigma}[\alpha^{2}(h_{22})\alpha(g_{2})]\alpha^{2}(l_{2})
	 \\=[a[\alpha^{m}(h_{1})\cdot(\beta^{-2}(b)(\alpha^{m}(g_{111})\cdot\beta^{-3}(c)))]][\sigma(\alpha^{k+3}(h_{211}),\alpha^{k+3}(g_{112}))\\~~~~~\sigma(\alpha^{k+2}(h_{212})\alpha^{k+1}(g_{12}),\alpha^{k+1}(l_{1}))]\sharp _{\sigma}[\alpha^{2}(h_{22})\alpha(g_{2})]\alpha^{2}(l_{2})
	 \\=[a[\alpha^{m}(h_{1})\cdot(\beta^{-2}(b)(\alpha^{m-2}(g_{1})\cdot\beta^{-3}(c)))]][\sigma(\alpha^{k+3}(h_{211}),\alpha^{k+3}(g_{211}))\\~~~~~\sigma(\alpha^{k+2}(h_{212}g_{212}),\alpha^{k+1}(l_{1}))]\sharp _{\sigma}\alpha^{2}(h_{22}g_{22})\alpha^{2}(l_{2})
	 \\\overset{(3)}{=}[a[\alpha^{m}(h_{1})\cdot(\beta^{-2}(b)(\alpha^{m-2}(g_{1})\cdot\beta^{-3}(c)))]][\alpha^{m+2}(h_{211})\cdot\sigma(\alpha^{k+2}(g_{211}),\alpha^{k+1}(l_{11}))]\\~~~~~\sigma(\alpha^{k+3}(h_{212}),\alpha^{k+2}(g_{212})\alpha^{k+1}(l_{12}))\sharp _{\sigma}\alpha^{2}(h_{22}g_{22})\alpha^{2}(l_{2})
	 \\=[[\beta^{-1}(a)[\alpha^{m+1}(h_{111})\cdot(\beta^{-3}(b)(\alpha^{m-3}(g_{1})\cdot\beta^{-4}(c)))]][\alpha^{m+2}(h_{112})\cdot\sigma(\alpha^{k+2}(g_{211}),\alpha^{k+1}(l_{11}))]]\\~~~~~\sigma(\alpha^{k+3}(h_{12}),\alpha^{k+3}(g_{212})\alpha^{k+2}(l_{12}))\sharp _{\sigma}[\alpha(h_{2})\alpha^{2}(g_{22})]\alpha^{2}(l_{2})
	 \\=[a[(\alpha^{m+1}(h_{111})\cdot(\beta^{-3}(b)(\alpha^{m-3}(g_{1})\cdot\beta^{-4}(c))))(\alpha^{m+1}(h_{112})\cdot\sigma(\alpha^{k+1}(g_{211}),\alpha^{k}(l_{11})))]]\\~~~~~\sigma(\alpha^{k+3}(h_{12}),\alpha^{k+3}(g_{212})^{k+2}(l_{12}))\sharp _{\sigma}[\alpha(h_{2})\alpha^{2}(g_{22})]\alpha^{2}(l_{2})
	 \\=[a[\alpha^{m+1}(h_{11})\cdot[[\beta^{-3}(b)(\alpha^{m-3}(g_{1})\cdot\beta^{-4}(c))]\sigma(\alpha^{k+1}(g_{211}),\alpha^{k}(l_{11}))]]]\\~~~~~\sigma(\alpha^{k+3}(h_{12}),\alpha^{k+3}(g_{212})\alpha^{k+2}(l_{12}))\sharp _{\sigma}[\alpha(h_{2})\alpha^{2}(g_{22})]\alpha^{2}(l_{2})
	 \\=\beta(a)[[\alpha^{m+1}(h_{11})\cdot[[\beta^{-3}(b)(\alpha^{m-2}(g_{11})\cdot\beta^{-4}(c))]\sigma(\alpha^{k}(g_{12}),\alpha^{k-1}(l_{1}))]]\\~~~~~\sigma(\alpha^{k+2}(h_{12}),\alpha^{k+1}(g_{21}l_{21}))]\sharp _{\sigma}[\alpha(h_{2})\alpha^{2}(g_{22})]\alpha^{3}(l_{22})
	 \\=\beta(a)[[\alpha^{m+1}(h_{11})\cdot[\beta^{-2}(b)[(\alpha^{m-2}(g_{11})\cdot\beta^{-4}(c))\sigma(\alpha^{k-1}(g_{12}),\alpha^{k-2}(l_{1}))]]]\\~~~~~\sigma(\alpha^{k+2}(h_{12}),\alpha^{k+1}(g_{21}l_{21}))]\sharp _{\sigma}\alpha^{2}(h_{2})\alpha^{2}(g_{22}l_{22})
	\\=[\beta(a)\sharp _{\sigma}\alpha(h)][b[(\alpha^{m}(g_{11})\cdot\beta^{-2}(c))\sigma(\alpha^{k+1}(g_{12}),\alpha^{k}(l_{1}))]\sharp _{\sigma}\alpha(g_{2}l_{2})]
	\\=[\beta(a)\sharp _{\sigma}\alpha(h)][(b\sharp _{\sigma} g)(c\sharp _{\sigma} l)] $.
	
	The proof is completed.
	$\hfill \blacksquare$
\\
	
{\bf Example 2.4.} Let $ H_{4}=\left\langle g, x | 1\cdot g=g, 1\cdot x=-x, g^{2}=1, x^{2}=0, xg=-gx \right\rangle $. Define the comultipication, counit and antipode by
	\begin{align*}
		\bigtriangleup(1)=1\otimes 1, ~~~&\varepsilon(1)=1, ~~~S(1)=1\\
		\bigtriangleup(g)=g\otimes g, ~~~&\varepsilon(g)=1, ~~~S(g)=g\\
		\bigtriangleup(x)=-x\otimes g+1\otimes(-x), ~~~&\varepsilon(x)=0, ~~~S(x)=-gx.\\
	\end{align*}
	Hom-map $ \alpha: H_{4}\rightarrow H_{4} $ is given by
	\begin{align*}
		\alpha(1)=1, \alpha(g)=g, \alpha(x)=-x, \alpha(gx)=-gx.
	\end{align*}
	Then $ (H_{4},\alpha) $ is a Hom-Hopf algebra.

	Let $ K[y] $ be a polynomial algebra over field $ K $ with the indeterminant $ y $ and $ K[y]/\left\langle y^{2}\right\rangle  $be the quotient algebra. Consider the Hom-algebra $ (K[y]/\left\langle y^{2}\right\rangle,id)  $ and define the action of $ H_{4} $ on  $ K[y]/\left\langle y^{2}\right\rangle  $ by
	\begin{align*}
		h\cdot 1=\varepsilon(h)1, 1\cdot y=y, g\cdot y=y, x\cdot y=0, gx\cdot y=0.
	\end{align*}
	Let $ n\in K $. Define a linear map $ \sigma:H_{4}\otimes H_{4}\rightarrow K[y]/\left\langle y^{2}\right\rangle $ by
\begin{table}[h]
\centering
\begin{tabular}{|c|c c c c|}
\hline
$ \sigma $ & 1 & $ g $ & $ x $ & $ gx $ \\
\hline
1&1&1&0&0\\
\hline
g&1&1&0&0\\
\hline
$ x $&0&0&$\frac{n}{2}$&-$\frac{n}{2}$\\
\hline
$ gx $&0&0&$\frac{n}{2}$&-$\frac{n}{2}$\\
\hline
\end{tabular}
\end{table}

Easy to see that $ \sigma $ satisfies the conditions in Theorem 2.3. Thus we have a $ (m,k) $-Hom-crossed product algebra $ (K[y]/\left\langle y^{2}\right\rangle\sharp_{\sigma}H_{4}, id\otimes \alpha) $.
\\

\section*{3. Radford $[(m,k),m]$-biproducts}
	\def\theequation{3. \arabic{equation}}
	\setcounter{equation} {0} \hskip\parindent
In this section, we will prove our Radford $[(m,k),m]$-biproduct theorem.

 Let $ (H,\alpha) $ be a monoidal Hom-bialgebra and $ (A,\beta) $ a left $ (H,\alpha) $-Hom-comodule coalgebra For any $ m\in \mathbb{Z }$.
 A left $m$-Hom-smash coproduct $( A\rtimes H,\beta\otimes\alpha) $ is  $ A\rtimes H=A\otimes H  $ as a linear space,
  with a Hom-comultiplication and counit
  given by:  $ \forall a\rtimes h\in A\rtimes H $,
	\begin{eqnarray*}
&&\v _{A\rtimes H}(a\rtimes h)=\v (a)\v (h),\\
	&& \bigtriangleup_{A\rtimes H}(a\rtimes h)=a_{1}\rtimes\alpha^{m}(a_{2(-1)})\alpha^{-1}(h_{1})\otimes\beta(a_{2(0)})\rtimes h_{2}.
	\end{eqnarray*}

{\bf Proposition 3.1.} With notations above. Then $( A\rtimes H,\beta\otimes\alpha) $ is a monoidal Hom-coalgebra.
	
	{\bf Proof.}  For all $ a\rtimes h\in A\rtimes H  $, it is easy to verify the following equation:
	
	\begin{center}
		$ \bigtriangleup_{A\rtimes H}(\beta(a)\rtimes \alpha(h))=(\beta\otimes\alpha)\bigtriangleup_{A\rtimes H}(a\rtimes h), $
		
		$ \varepsilon_{A\rtimes H}(\beta(a)\rtimes\alpha(h))=\varepsilon_{A}(a)\varepsilon_{H}(h)=\varepsilon_{A\rtimes H}(a\rtimes h), $
		
	 $(id\otimes\varepsilon_{A\rtimes H})\bigtriangleup_{A\rtimes H}(a\rtimes h)=\beta^{-1}(a)\otimes\alpha^{-1}(h)=(\varepsilon_{A\rtimes H}\otimes id)\bigtriangleup_{A\rtimes H}(a\rtimes h). $
	\end{center}
	
	Besides, we also have
	
	$ (\beta^{-1}\otimes\alpha^{-1})(a_{1}\rtimes\alpha^{m}(a_{2(-1)})\alpha^{-1}(h_{1}))\otimes\bigtriangleup_{A\rtimes H}(\beta(a_{2(0)})\rtimes h_{2})
	 \\=\beta^{-1}(a_{1})\rtimes\alpha^{m-1}(a_{2(-1)})\alpha^{-2}(h_{1})\otimes\beta(a_{2(0)1})\rtimes\alpha^{m+1}(a_{2(0)2(-1)})\alpha^{-1}(h_{21})\\~~~~~\otimes\beta^{2}(a_{2(0)2(0)})\rtimes h_{22}
	 \\=\beta^{-1}(a_{1})\rtimes\alpha^{m-1}(a_{2(-1)})\alpha^{-1}(h_{11})\otimes\beta(a_{2(0)1})\rtimes\alpha^{m+1}(a_{2(0)2(-1)})\alpha^{-1}(h_{12})\\~~~~~\otimes\beta^{2}(a_{2(0)2(0)})\rtimes \alpha^{-1}(h_{2})
	 \\=\beta^{-1}(a_{1})\rtimes\alpha^{m-1}(a_{21(-1)}a_{22(-1)})\alpha^{-1}(h_{11})\otimes\beta(a_{21(0)})\rtimes\alpha^{m+1}(a_{22(0)(-1)})\alpha^{-1}(h_{12})\\~~~~~\otimes\beta^{2}(a_{22(0)(0)})\rtimes \alpha^{-1}(h_{2})
	 \\=\beta^{-1}(a_{1})\rtimes[\alpha^{m-1}(a_{21(-1)})\alpha^{m}(a_{22(-1)1})]\alpha^{-1}(h_{11})\otimes\beta(a_{21(0)})\rtimes\alpha^{m+1}(a_{22(-1)2})\alpha^{-1}(h_{12})\\~~~~~\otimes\beta^{2}(a_{22(0)})\rtimes \alpha^{-1}(h_{2}) $.
	
	$ \bigtriangleup_{A\rtimes H}(a_{1}\rtimes\alpha^{m}(a_{2(-1)})\alpha^{-1}(h_{1}))\otimes(\beta^{-1}\otimes\alpha^{-1})(\beta(a_{2(0)})\rtimes h_{2})
	 \\=a_{11}\rtimes\alpha^{m}(a_{12(-1)})[\alpha^{m-1}(a_{2(-1)1})\alpha^{-2}(h_{11})]\otimes\beta(a_{12(0)})\rtimes\alpha^{m}(a_{2(-1)2})\alpha^{-1}(h_{12})\\~~~~~\otimes a_{2(0)}\rtimes\alpha^{-1}(h_{2})
	 \\=\beta^{-1}(a_{1})\rtimes\alpha^{m}(a_{21(-1)})[\alpha^{m}(a_{22(-1)1})\alpha^{-2}(h_{11})]\otimes\beta(a_{21(0)})\rtimes\alpha^{m+1}(a_{22(-1)2})\alpha^{-1}(h_{12})\\~~~~~\otimes a_{22(0)}\rtimes\alpha^{-1}(h_{2})
	 \\=\beta^{-1}(a_{1})\rtimes[\alpha^{m-1}(a_{21(-1)})\alpha^{m}(a_{22(-1)1})]\alpha^{-1}(h_{11})\otimes\beta(a_{21(0)})\rtimes\alpha^{m+1}(a_{22(-1)2})\alpha^{-1}(h_{12})\\~~~~~\otimes\beta^{2}(a_{22(0)})\rtimes \alpha^{-1}(h_{2}) $
	
Thus,
	\begin{center}
		$  ((\beta^{-1}\otimes\alpha^{-1})\otimes\bigtriangleup_{A\rtimes H})\circ\bigtriangleup_{A\rtimes H}=(\bigtriangleup_{A\rtimes H}\otimes(\beta^{-1}\otimes\alpha^{-1}))\circ\bigtriangleup_{A\rtimes H} $.
	\end{center}

	The proof is completed.
		$\hfill \blacksquare$
\\
	
	{\bf Definition 3.2.}  Let $ (A\sharp_{\sigma} H,\beta\otimes\alpha) $ be $ (m,k) $-Hom-crossed product.
 Then $ \sigma $ is called a {\sl twisted comodule cocycle} if
	\begin{align*}
		\beta(a_{1})\otimes\alpha^{m+1}(a_{2(-1)})g\otimes a_{2(0)}=a_{1}\sigma(\alpha^{k+m+2}(a_{2(-1)1}),\alpha^{k+1}(g_{1}))\otimes\alpha^{m+2}(a_{2(-1)2})\alpha(g_{2})\otimes a_{2(0)}.
	\end{align*}

   {\bf Theorem 3.3.}  Let $ (H,\alpha) $ be a monoidal Hom-bialgebra and $ (A,\beta) $ a monoidal Hom-algebra. Suppose that $ (H,\alpha) $ weakly acts on $ (A,\beta) $ and $ (A,\beta) $ is a left $ (H,\alpha) $ Hom-comodule coalgebra with the comodule structure map $ \rho_{A}^{H}:A\rightarrow H\otimes A$. If $ (A\sharp_{\sigma} H,\beta\otimes\alpha) $ is a $ (m,k) $-Hom-crossed product with $ \sigma $ being a twisted comodule cocycle and
   $(A\rtimes H,\beta\otimes\alpha) $ is a $ m $-Hom-smash coproduct, then the following conditions are equivalent:

   (1) $ (A^{\sharp_{\sigma}}_{\rtimes} H,m_{A\sharp_{\sigma} H},1_{A}\otimes 1_{H},\bigtriangleup_{A\rtimes H},\varepsilon_{A\rtimes H},\beta\otimes\alpha) $
     is a monoidal  Hom-bialgebra.

   (2)  The conditions:

   \quad (A1) $ \varepsilon_{A} $ is a Hom-algebra map,

   \quad (A2) $ \varepsilon_{A}(h\cdot a)=\varepsilon_{H}(h)\varepsilon_{A}(a) $,

   \quad (A3) $ \sigma $ is a Hom-coalgebra map,

   \quad (A4) $ \bigtriangleup_{A}(1_{A})=1_{A}\otimes 1_{A} $,

   \quad (A5) $ \rho_{A}^{H}(ab)=\rho_{A}^{H}(a)\rho_{A}^{H}(b) , $
   $ \rho_{A}^{H}(1_{A})=1_{H}\otimes 1_{A} $,

   \quad (A6) $ \alpha^{m-1}(\sigma(\alpha^{k+2}(h_{1}),\alpha^{k+2}(g_{1}))_{(-1)})\alpha^{-1}(h_{2}g_{2})\otimes\sigma(\alpha^{k+2}(h_{1}),\alpha^{k+2}(g_{1}))_{(0)}
   =h_{1}g_{1}\otimes\sigma(\alpha^{k+1}(h_{2}),\alpha^{k+1}(g_{2})) $,

  \quad (A7) $ \bigtriangleup_{A}(ab)=a_{1}[(\alpha^{2m}(a_{2(-1)1})\cdot\beta^{-2}(b_{1}))\sigma(\alpha^{k+m+1}(a_{2(-1)2}),\alpha^{k+m}(b_{2(-1)}))]\otimes\beta(a_{2(0)}b_{2(0)}) $,

  \quad (A8)
   $ \bigtriangleup_{A}(\alpha^{m}(h)\cdot b)=(\alpha^{m}(h_{11})\cdot\beta^{-1}(b_{1}))\sigma(\alpha^{k+1}(h_{12}),\alpha^{k+m+1}(b_{2(-1)}))\otimes\alpha^{m}(h_{2})\cdot\beta(b_{2(0)}) $,

  \quad (A9)
   $ \alpha^{m-1}[(\alpha^{m+1}(h_{1})\cdot b)_{(-1)}]h_{2}\otimes(\alpha^{m+1}(h_{1})\cdot b)_{(0)}=h_{1}\alpha^{m}(b_{(-1)})\otimes\alpha^{m}(h_{2})\cdot b_{(0)}$,\\
    for any $a, b\in A$ and $h\in H$.

{\bf Proof.} (1)$\Rightarrow$(2) follows from the similar calculations to those of
 \cite[Theorem 1]{WJZ}. So we only need to show (2)$\Rightarrow$(1). Assume (2) holds, then $ \varepsilon_{A\rtimes H} $ is Hom-algebra map by (A1) and (A2). By (A4) and (A5), we have $ \bigtriangleup_{A\rtimes H}(1_{A}\sharp_{\sigma}1_{H})=1_{A}\sharp_{\sigma}1_{H}\otimes1_{A}\sharp_{\sigma}1_{H} $.

   Next we prove:
   \begin{align*}
   	\bigtriangleup_{A\rtimes H}[(a\sharp_{\sigma} h)(b\sharp_{\sigma} g)]=\bigtriangleup_{A\rtimes H}(a\sharp_{\sigma} h)\bigtriangleup_{A\rtimes H}(b\sharp_{\sigma} g).
   \end{align*}
   It is enough to verify the following relations:
   \begin{align}
   	\bigtriangleup_{A\rtimes H}[(a\sharp_{\sigma} 1_{H})(b\sharp_{\sigma} 1_{H})]=\bigtriangleup_{A\rtimes H}(a\sharp_{\sigma} 1_{H})\bigtriangleup_{A\rtimes H}(b\sharp_{\sigma} 1_{H}) ;
   	\\  \bigtriangleup_{A\rtimes H}[(a\sharp_{\sigma} 1_{H})(1_{A}\sharp_{\sigma} g)]=\bigtriangleup_{A\rtimes H}(a\sharp_{\sigma} 1_{H})\bigtriangleup_{A\rtimes H}(1_{A}\sharp_{\sigma} g);
   	\\  \bigtriangleup_{A\rtimes H}[(1_{A}\sharp_{\sigma} h)(b\sharp_{\sigma} 1_{H})]=\bigtriangleup_{A\rtimes H}(1_{A}\sharp_{\sigma} h)\bigtriangleup_{A\rtimes H}(b\sharp_{\sigma} 1_{H}) ;
   	\\  \bigtriangleup_{A\rtimes H}[(1_{A}\sharp_{\sigma} h)(1_{A}\sharp_{\sigma} g)]=\bigtriangleup_{A\rtimes H}(1_{A}\sharp_{\sigma} h)\bigtriangleup_{A\rtimes H}(1_{A}\sharp_{\sigma} g) .
   \end{align}
   In fact, if (3.1)-(3.4) holds, then
   \begin{align*}
   	\bigtriangleup_{A\rtimes H}[(a\sharp_{\sigma} 1_{H})(b\sharp_{\sigma} g)]
   	=&\bigtriangleup_{A\rtimes H}[ab\sharp_{\sigma}\alpha(g)]\\=&\bigtriangleup_{A\rtimes H}[(\beta^{-1}(ab)\sharp_{\sigma} 1_{H})(1_{A}\sharp_{\sigma} g)]
   	\\\overset{(3.2)}{=}&\bigtriangleup_{A\rtimes H}(\beta^{-1}(ab)\sharp_{\sigma} 1_{H})\bigtriangleup_{A\rtimes H}(1_{A}\sharp_{\sigma} g)
   	\\=&\bigtriangleup_{A\rtimes H}[(\beta^{-1}(a)\sharp_{\sigma} 1_{H})(\beta^{-1}(b)\sharp_{\sigma} 1_{H})]\bigtriangleup_{A\rtimes H}(1_{A}\sharp_{\sigma} g)
   	\\\overset{(3.1)}{=}&[\bigtriangleup_{A\rtimes H}(\beta^{-1}(a)\sharp_{\sigma} 1_{H})\bigtriangleup_{A\rtimes H}(\beta^{-1}(b)\sharp_{\sigma} 1_{H})]\bigtriangleup_{A\rtimes H}(1_{A}\sharp_{\sigma} g)
   	\\=&\bigtriangleup_{A\rtimes H}(a\sharp_{\sigma} 1_{H})[\bigtriangleup_{A\rtimes H}(\beta^{-1}(b)\sharp_{\sigma} 1_{H})\bigtriangleup_{A\rtimes H}(1_{A}\sharp_{\sigma} \alpha^{-1}(g))]
   	\\\overset{(3.2)}{=}&\bigtriangleup_{A\rtimes H}(a\sharp_{\sigma} 1_{H})\bigtriangleup_{A\rtimes H}[(\beta(b)\sharp_{\sigma} 1_{H})(1_{A}\sharp_{\sigma} \alpha^{-1}(g))] \\=&\bigtriangleup_{A\rtimes H}(a\sharp_{\sigma} 1_{H})\bigtriangleup_{A\rtimes H}(b\sharp_{\sigma} g)
   \end{align*}
   and
   \begin{align*}
   	\bigtriangleup_{A\rtimes H}[(a\sharp_{\sigma} h)(1_{A}\sharp_{\sigma} g)]
   	=&\bigtriangleup_{A\rtimes H}[[(\beta^{-1}(a)\sharp_{\sigma} 1_{H})(1_{A}\sharp_{\sigma}\alpha^{-1}(h))](1_{A}\sharp_{\sigma} g)]
   	\\=&\bigtriangleup_{A\rtimes H}[(a\sharp_{\sigma} 1_{H})[(1_{A}\sharp_{\sigma}\alpha^{-1}(h))(1_{A}\sharp_{\sigma}\alpha^{-1}(g))]]
   	\\=&\bigtriangleup_{A\rtimes H}(a\sharp_{\sigma} 1_{H})\bigtriangleup_{A\rtimes H}[(1_{A}\sharp_{\sigma}\alpha^{-1}(h))(1_{A}\sharp_{\sigma}\alpha^{-1}(g))]
   	\\\overset{(3.4)}{=}&\bigtriangleup_{A\rtimes H}(a\sharp_{\sigma} 1_{H})[\bigtriangleup_{A\rtimes H}(1_{A}\sharp_{\sigma}\alpha^{-1}(h))\bigtriangleup_{A\rtimes H}(1_{A}\sharp_{\sigma}\alpha^{-1}(g))]
   	\\=&[\bigtriangleup_{A\rtimes H}(\beta^{-1}(a)\sharp_{\sigma} 1_{H})\bigtriangleup_{A\rtimes H}(1_{A}\sharp_{\sigma}\alpha^{-1}(h))]\bigtriangleup_{A\rtimes H}(1_{A}\sharp_{\sigma} g)
   	\\\overset{(3.2)}{=}&\bigtriangleup_{A\rtimes H}[(\beta^{-1}(a)\sharp_{\sigma} 1_{H})(1_{A}\sharp_{\sigma}\alpha^{-1}(h))]\bigtriangleup_{A\rtimes H}(1_{A}\sharp_{\sigma} g)
   	\\=&\bigtriangleup_{A\rtimes H}(a\sharp_{\sigma} h)\bigtriangleup_{A\rtimes H}(1_{A}\sharp_{\sigma} g).
   \end{align*}

Then we have
\begin{align*}
	\bigtriangleup_{A\rtimes H}[(a\sharp_{\sigma} h)(b\sharp_{\sigma} g)]
	=&\bigtriangleup_{A\rtimes H}[(a\sharp_{\sigma} h)[(\beta^{-1}(b)\sharp_{\sigma} 1_{H})(1_{A}\sharp_{\sigma}\alpha^{-1}(g))]]
	\\=&\bigtriangleup_{A\rtimes H}[[(\beta^{-1}(a)\sharp_{\sigma}\alpha^{-1}(h))(\beta^{-1}(b)\sharp_{\sigma} 1_{H})](1_{A}\sharp_{\sigma} g)]
	\\=&\bigtriangleup_{A\rtimes H}[(\beta^{-1}(a)\sharp_{\sigma}\alpha^{-1}(h))(\beta^{-1}(b)\sharp_{\sigma} 1_{H})]\bigtriangleup_{A\rtimes H}(1_{A}\sharp_{\sigma} g)
	\\=&\bigtriangleup_{A\rtimes H}[[(\beta^{-2}(a)\sharp_{\sigma} 1_{H})(1_{A}\sharp_{\sigma}\alpha^{-2}(h))](\beta^{-1}(b)\sharp_{\sigma} 1_{H})]\bigtriangleup_{A\rtimes H}(1_{A}\sharp_{\sigma} g)
	\\=&[\bigtriangleup_{A\rtimes H}(\beta^{-1}(a)\sharp_{\sigma} 1_{H})\bigtriangleup_{A\rtimes H}[(1_{A}\sharp_{\sigma}\alpha^{-2}(h))(\beta^{-2}(b)\sharp_{\sigma} 1_{H})]]\bigtriangleup_{A\rtimes H}(1_{A}\sharp_{\sigma} g)
	\\\overset{(3.3)}{=}&[(\beta^{-1}(a)\sharp_{\sigma} 1_{H})[\bigtriangleup_{A\rtimes H}(1_{A}\sharp_{\sigma}\alpha^{-2}(h))\bigtriangleup_{A\rtimes H}(\beta^{-2}(b)\sharp_{\sigma} 1_{H})]]\bigtriangleup_{A\rtimes H}(1_{A}\sharp_{\sigma} g)
	\\\overset{(3.2)}{=}&[\bigtriangleup_{A\rtimes H}[(\beta^{-1}(a)\sharp_{\sigma}\alpha^{-1}(h))]\bigtriangleup_{A\rtimes H}(\beta^{-1}(b)\sharp_{\sigma} 1_{H})]\bigtriangleup_{A\rtimes H}(1_{A}\sharp_{\sigma} g)
	\\=&\bigtriangleup_{A\rtimes H}(a\sharp_{\sigma} h)[\bigtriangleup_{A\rtimes H}(\beta^{-1}(b)\sharp_{\sigma} 1_{H})\bigtriangleup_{A\rtimes H}(1_{A}\sharp_{\sigma}\alpha^{-1}(g))]
	\\\overset{(3.2)}{=}&\bigtriangleup_{A\rtimes H}(a\sharp_{\sigma} h)\bigtriangleup_{A\rtimes H}(b\sharp_{\sigma} g).
\end{align*}

Next we show that (2.1)-(2.4) hold, that is
\begin{align*}
	&\bigtriangleup_{A\rtimes H}(a\sharp_{\sigma} 1_{H})\bigtriangleup_{A\rtimes H}(b\sharp_{\sigma} 1_{H})
	 \\=&a_{1}[(\alpha^{2m+1}(a_{2(-1)11})\cdot\beta^{-2}(b_{1}))\sigma(\alpha^{k+m+2}(a_{2(-1)12}),\alpha^{k+m+1}(b_{2(-1)1}))]\\&\sharp_{\sigma}\alpha^{m+2}(a_{2(-1)2}b_{2(-1)2})\otimes\beta(a_{2(0)}b_{2(0)})\sharp_{\sigma} 1_{H}
	 \\=&a_{1}[(\alpha^{2m}(a_{2(-1)1})\cdot\beta^{-2}(b_{1}))\sigma(\alpha^{k+m+1}(a_{2(-1)2}),\alpha^{k+m}(b_{2(-1)}))]\\&\sharp_{\sigma}\alpha^{m+2}(a_{2(0)(-1)}b_{2(0)(-1)})\otimes\beta(a_{2(0)(0)}b_{2(0)(0)})\sharp_{\sigma} 1_{H}
	 \\\overset{(A7)}{=}&(ab)_{1}\sharp_{\sigma}\alpha^{m+1}[(ab)_{2(-1)}]\otimes\beta(ab)_{2(0)}\sharp_{\sigma} 1_{H}
	\\=&\bigtriangleup_{A\rtimes H}(ab\sharp_{\sigma} 1_{H}),
\end{align*}
and (3.1) is proved.
\begin{align*}
	&\bigtriangleup_{A\rtimes H}(a\sharp_{\sigma} 1_{H})(1_{A}\sharp_{\sigma} g)
	 \\=&a_{1}\sigma(\alpha^{k+m+2}(a_{2(-1)1}),\alpha^{k+1}(g_{11}))\sharp_{\sigma}\alpha^{m+2}(a_{2(-1)2})\alpha(g_{12})\otimes\beta^{2}(a_{2(0)})\sharp_{\sigma}\alpha(g_{2})
	 \\=&\beta(a_{1})\sharp_{\sigma}\alpha^{m+1}(a_{2(-1)})g_{1}\otimes\beta^{2}(a_{2(0)})\sharp_{\sigma}\alpha(g_{2})
	\\=&\bigtriangleup_{A\rtimes H}(\beta(a)\sharp_{\sigma}\alpha(g))
	\\=&\bigtriangleup_{A\rtimes H}[(a\sharp_{\sigma} 1_{H})(1_{A}\sharp_{\sigma} g)],
\end{align*}
	and (3.2) is proved.
	\begin{align*}
		&\bigtriangleup_{A\rtimes H}(1_{A}\sharp_{\sigma} h)\bigtriangleup_{A\rtimes H}(b\sharp_{\sigma} 1_{H})
		\\= &(\alpha^{m+1}(h_{111})\cdot\beta^{-1}(b_{1}))\sigma(\alpha^{k+2}(h_{112}),\alpha^{k+m+2}(b_{2(-1)1}))\sharp_{\sigma}\alpha(h_{12})\alpha^{m+2}(b_{2(-1)2})\\&\otimes(\alpha^{m+1}(h_{21})\cdot\beta(b_{2(0)}))\sharp_{\sigma}\alpha^{2}(h_{22})
		 \\=&(\alpha^{m}(h_{11})\cdot\beta^{-1}(b_{1}))\sigma(\alpha^{k+1}(h_{12}),\alpha^{k+m+2}(b_{2(-1)1}))\sharp_{\sigma}\alpha(h_{21})\alpha^{m+2}(b_{2(-1)2})\\&\otimes(\alpha^{m+2}(h_{221})\cdot\beta(b_{2(0)}))\sharp_{\sigma}\alpha^{3}(h_{222})
		 \\=&(\alpha^{m}(h_{11})\cdot\beta^{-1}(b_{1}))\sigma(\alpha^{k+1}(h_{12}),\alpha^{k+m+2}(b_{2(-1)1}))\sharp_{\sigma}\alpha^{2}(h_{211})\alpha^{m+2}(b_{2(-1)2})\\&\otimes(\alpha^{m+2}(h_{212})\cdot\beta(b_{2(0)}))\sharp_{\sigma}\alpha^{2}(h_{22})
		 \\=&(\alpha^{m}(h_{11})\cdot\beta^{-1}(b_{1}))\sigma(\alpha^{k+1}(h_{12}),\alpha^{k+m+1}(b_{2(-1)}))\sharp_{\sigma}\alpha^{2}(h_{211})\alpha^{m+2}(b_{2(0)(-1)})\\&\otimes(\alpha^{m+2}(h_{212})\cdot\beta^{2}(b_{2(0)(0)}))\sharp_{\sigma}\alpha^{2}(h_{22})
		 \\\overset{(A9)}{=}&(\alpha^{m}(h_{11})\cdot\beta^{-1}(b_{1}))\sigma(\alpha^{k+1}(h_{12}),\alpha^{k+m+1}(b_{2(-1)}))\sharp_{\sigma}\alpha[(\alpha^{m+2}(h_{211})\cdot\beta(b_{2(0)}))_{(-1)}]\\&\alpha^{2}(h_{212})\otimes\beta[(\alpha^{m+2}(h_{211})\cdot\beta(b_{2(0)}))_{(0)}]\sharp_{\sigma}\alpha^{2}(h_{22})
		 \\=&(\alpha^{m}(h_{11})\cdot\beta^{-1}(b_{1}))\sigma(\alpha^{k+1}(h_{12}),\alpha^{k+m+1}(b_{2(-1)}))\sharp_{\sigma}\alpha^{m}[(\alpha^{m+2}(h_{211})\cdot\beta(b_{2(0)}))_{(-1)}]\\&\alpha^{2}(h_{212})\otimes\beta[(\alpha^{m+2}(h_{211})\cdot\beta(b_{2(0)}))_{(0)}]\sharp_{\sigma}\alpha^{2}(h_{22})
		 \\=&(\alpha^{m+1}(h_{111})\cdot\beta^{-1}(b_{1}))\sigma(\alpha^{k+2}(h_{112}),\alpha^{k+m+1}(b_{2(-1)}))\sharp_{\sigma}\alpha^{m}[(\alpha^{m+1}(h_{12})\cdot\beta(b_{2(0)}))_{(-1)}]\\&\alpha(h_{21})\otimes\beta[(\alpha^{m+1}(h_{12})\cdot\beta(b_{2(0)}))_{(0)}]\sharp_{\sigma}\alpha^{2}(h_{22})
		\\\overset{(A8)}{=}&(\alpha^{m+1}(h_{1})\cdot b)_{1}\sharp_{\sigma}\alpha^{m}[(\alpha^{m+1}(h_{1})\cdot b)_{2(-1)}]\alpha(h_{21})\otimes\beta[(\alpha^{m+1}(h_{1})\cdot b)_{2(0)}]\sharp_{\sigma}\alpha^{2}(h_{22})
		\\=&\bigtriangleup_{A\rtimes H}[(\alpha^{m+1}(h_{1})\cdot b)\sharp_{\sigma}\alpha^{2}(h_{2})]
		\\=&\bigtriangleup_{A\rtimes H}[(1_{A}\sharp_{\sigma} h)(b\sharp_{\sigma} 1_{H})],
	\end{align*}
	and (3.3) is proved.
	\begin{align*}
		&\bigtriangleup_{A\rtimes H}(1_{A}\sharp_{\sigma} h)\bigtriangleup_{A\rtimes H}(1_{A}\sharp_{\sigma} g)
		 \\=&\sigma(\alpha^{k+2}(h_{11}),\alpha^{k+2}(g_{11}))\sharp_{\sigma}\alpha(h_{12}g_{12})\otimes\sigma(\alpha^{k+2}(h_{21}),\alpha^{k+2}(g_{21}))\sharp_{\sigma}\alpha(h_{22}g_{22})
		 \\=&\sigma(\alpha^{k+1}(h_{1}),\alpha^{k+1}(g_{1}))\sharp_{\sigma}\alpha(h_{211}g_{211})\otimes\sigma(\alpha^{k+3}(h_{212}),\alpha^{k+3}(g_{212}))\sharp_{\sigma}\alpha(h_{22}g_{22})
		 \\\overset{(A6)}{=}&\sigma(\alpha^{k+2}(h_{11}),\alpha^{k+2}(g_{11}))\sharp_{\sigma}\alpha^{m}(\sigma(\alpha^{k+2}(h_{1}),\alpha^{k+2}(g_{1}))_{2(-1)})(h_{21}g_{21})\\&\otimes\beta(\sigma(\alpha^{k+2}(h_{1}),\alpha^{k+2}(g_{1}))_{2(0)})\sharp_{\sigma}\alpha(h_{22}g_{22})
		\\=&\bigtriangleup_{A\rtimes H}[\sigma(\alpha^{k+2}(h_{1}),\alpha^{k+2}(g_{1}))\sharp_{\sigma}\alpha(h_{2}g_{2})]
		\\=&\bigtriangleup_{A\rtimes H}[(1_{A}\sharp_{\sigma} h)(1_{A}\sharp_{\sigma} g)].
	\end{align*}
	and (3.4) is proved.
	
	Therefore, $ (A^{\sharp_{\sigma}}_{\rtimes} H,m_{A\sharp_{\sigma} H},1_{A}\otimes 1_{H},\bigtriangleup_{A\rtimes H},\varepsilon_{A\rtimes H},\beta\otimes\alpha) $ is a monoidal Hom-bialgebra.
	$\hfill \blacksquare$
\\

{\bf Definition 3.4.}  Let $ (H,\alpha) $ be a monoidal Hom-bialgebra and $ (A,\beta) $ a monoidal Hom-algebra. Let $\sigma:H\otimes H\rightarrow A$ and $S: H\rightarrow H $ be linear maps. Then $S $ is called a $ \sigma $-antipode of $ (H,\alpha) $ if
	
	(1) $ \alpha\circ S=S\circ\alpha $,
	
	(2) $(\sigma\otimes m_{H})\bigtriangleup_{H\otimes H}(id_{H}\otimes S)\bigtriangleup_{H}(h)=\varepsilon_{H}(h)1_{A}\otimes1_{H} $,
	
	(3) $ (\sigma\otimes m_{H})\bigtriangleup_{H\otimes H}(S\otimes id_{H})\bigtriangleup_{H}(h)=\varepsilon_{H}(h)1_{A}\otimes1_{H} $.
	\\
In this case $ (H,\alpha) $ is called a $ \sigma $-monoidal Hom-Hopf algebra.
\\

{\bf Theorem 3.5.} Let $ (A\sharp_{\sigma} H,\beta\otimes\alpha) $ be a monoidal  Hom-bialgebra. If $ (H,\alpha) $ is a $ \sigma $-monoidal Hom-Hopf algebra with $ \sigma $-antipode $ S_{H} $ and $ S_{A}\in Hom(A,A) $ is a convolution invertible element of $id_{A} $ with $ \beta\circ S_{A}=S_{A}\circ\beta $.
 Then $ (A^{\sharp_{\sigma}}_{\rtimes} H,m_{A\sharp_{\sigma} H},1_{A}\otimes 1_{H},\bigtriangleup_{A\rtimes H},\varepsilon_{A\rtimes H},\beta\otimes\alpha) $ is a monoidal Hom-Hopf algebra with the antipode given by
	\begin{align*}
		S_{A^{\sharp_{\sigma}}_{\rtimes} H}(a\otimes h)=(1_{A}~\otimes S_{H}(\alpha^{m-1}(a_{(-1)})\alpha^{-2}(h))(S_{A}(a_{(0)})\otimes 1_{H}).
	\end{align*}
	
	{\bf Proof.} For any $ a\otimes h\in A^{\sharp_{\sigma}}_{\rtimes} H  $, it is easy to verify the following equation:
	
	\begin{center}
		$ S_{A^{\sharp_{\sigma}}_{\rtimes} H}(\beta(a)\otimes \alpha(h))=(\beta\otimes\alpha)S_{A^{\sharp_{\sigma}}_{\rtimes} H}(a\otimes h). $
	\end{center}

Furthermore, we have
\begin{eqnarray*}	
&& S_{A^{\sharp_{\sigma}}_{\rtimes} H}(a_{1}~\otimes \alpha^{m}(a_{2(-1)})\alpha^{-1}(h_{1}))(\beta(a_{2(0)})\otimes h_{2})
	\\
&=&[[1_{A}~\otimes S_{H}[\alpha^{m-1}(a_{1(-1)})\alpha^{-2}(\alpha^{m}(a_{2(-1)})\alpha^{-1}(h_{1}))]](S_{A}(a_{1(0)})\otimes 1_{H})](\beta(a_{2(0)})\otimes h_{2})
\\
&=&[1_{A}~\otimes S_{H}[\alpha^{m}(a_{1(-1)})\alpha^{-1}(\alpha^{m}(a_{2(-1)})\alpha^{-1}(h_{1}))]][(S_{A}(a_{1(0)})\otimes 1_{H})(a_{2(0)}~\otimes \alpha^{-1}(h_{2}))]
	\\
&=&[1_{A}~\otimes S_{H}[\alpha^{m}(a_{1(-1)})\alpha^{-1}(\alpha^{m}(a_{2(-1)})\alpha^{-1}(h_{1}))]](S_{A}(a_{1(0)})a_{2(0)}~\otimes h_{2})
	\\
&=&[1_{A}~\otimes S_{H}[\alpha^{m-1}(a_{1(-1)}a_{2(-1)})\alpha^{-1}(h_{1})]](S_{A}(a_{1(0)})a_{2(0)}~\otimes h_{2})
	\\
&=&[1_{A}~\otimes S_{H}[\alpha^{m-1}(a_{(-1)})\alpha^{-1}(h_{1})]](S_{A}(a_{(0)1})a_{(0)2}~\otimes h_{2})
	\\
&=&\varepsilon_{A}(a)(1_{A}~\otimes S_{H}(h_{1}))(1_{A}~\otimes h_{2})
	 \\
&=&\varepsilon_{A}(a)\beta^{k+2}\circ\sigma(S(h_{1})_{1},h_{21})\otimes\alpha[S(h_{1})_{2}h_{22}]
	\\
&=&\varepsilon_{A}(a)\varepsilon_{H}(h)1_{A}~\otimes 1_{H}.
\end{eqnarray*}
Similarly, we also have
$$
(a_{1}\otimes\alpha^{m}(a_{2(-1)})\alpha^{-1}(h_{1}))S_{A^{\sharp_{\sigma}}_{\rtimes} H}(\beta(a_{2(0)})\otimes h_{2})=\varepsilon_{A}(a)\varepsilon_{H}(h)1_{A}\otimes 1_{H}.
$$
	
	 Therefore, $ S_{A^{\sharp_{\sigma}}_{\rtimes} H} $ is the convolution inverse of $ id_{A^{\sharp_{\sigma}}_{\rtimes} H} $ and $ (A^{\sharp_{\sigma}}_{\rtimes} H,m_{A\sharp_{\sigma} H},1_{A}\otimes 1_{H},\bigtriangleup_{A\rtimes H},\varepsilon_{A\rtimes H},\beta\otimes\alpha) $ is a monoidal Hom-Hopf algebra.
	 $\hfill \blacksquare$

{\bf Remark 3.6.} In the case of Hopf algebras, it follows from Theorem 3.5 that \cite[Theorem 2.5]{WJZ} by
  taking $m=0$ and  $k=-1$,  and  furthermore, it follows that \cite[Theorem 1]{R85} when $\sigma$ is trivial.

\section*{4. Hom admissible mappping system}
	\def\theequation{4. \arabic{equation}}
	\setcounter{equation}{0} \hskip\parindent

Let $ (C_{\sharp_{\sigma}}^{\times}H, \b \o \a ):=(C^{\sharp_{\sigma}}_{\rtimes} H, m_{A\sharp_{\sigma} H}, 1_{A}\otimes 1_{H},
 \bigtriangleup_{C\rtimes H},\varepsilon_{C\rtimes H},\beta\otimes\alpha)$.  In this section, we characterize
 the structure of the monoidal Hom-Hopf algebra $(C_{\sharp_{\sigma}}^{\times}H, \beta\otimes\alpha)$ established in Theorem 3.3.

We first claim that $\sigma$ is always convolution invertible with inverse $\sigma^{-1}$ in this section. It is not hard  to verify  the following Lemma.

\textbf{Lemma 4.1.} Let $(C\sharp_{\sigma}H,\beta\otimes\alpha)$ be a $(m,k)$ Hom-crossed product with convolution invertible $\sigma$.
Then

(1) $\alpha(h)\rightarrow\sigma(\alpha^{k+1}(l),\alpha^{k+1}(g))=[\sigma(\alpha^{k+2}(h_{11}),\alpha^{k+2}(l_{11}))\sigma(\alpha^{k+1}(h_{12}l_{12}),
\alpha^{k+1}(g_{1}))]\\\sigma^{-1}(\alpha^{k+2}(h_{2}),\alpha^{k+1}(l_{2}g_{2}))$,

(2) $\alpha^{m}(hl)\rightarrow\beta^{2}(a)=[\sigma(\alpha^{k+2}(h_{11}),\alpha^{k+2}(l_{11}))(\alpha^{m}(h_{12}l_{12})\rightarrow a)]\sigma^{-1}(\alpha^{k+2}(h_{2}),\alpha^{k+2}(l_{2}))$.
\\

\textbf{Definition 4.2.} Let $(H,\alpha)$ be a monoidal Hom Hopf algebra and $(C,\beta)$ a monoidal Hom coalgebra.
 Suppose $\bar{\sigma}: H\otimes H\lr C$ is bilinear and $\bar{\rho}: H\otimes H\lr C\otimes H\otimes H$ is a bilinear map defined by
\begin{equation*}
  \bar{\rho}(h\otimes l)=\bar{\sigma}(h_{1},l_{1})\otimes h_{2}\otimes l_{2}.
\end{equation*}
If there is a map $\varphi^{l}: H\otimes C\lr C$ such that

(1) $\varphi^{l}(\alpha\otimes\varphi^{l})=m_{C}(\beta\otimes\varphi^{l})(id\otimes m_{H}\otimes id)(\bar{\rho}\otimes id)$,

(2) $\varphi^{l}(1_{H}\otimes c)=\beta(c)$
\\
hold, then $(C,\beta)$ is called left $(H,\alpha,\bar{\sigma})$-Hom module.

\textbf{Remark 4.3.} If there is a map $\varphi^{r}: C\otimes H\rightarrow C$ such that

(1) $\varphi^{r}(\varphi^{r}\otimes\alpha)=m_{C}(\beta\otimes\varphi^{r})(id\otimes id\otimes m_{H})(id\otimes\bar{\rho})$,

(2) $\varphi^{r}(c\otimes 1_{H} )=\beta(c)$
\\
hold, then $(C,\beta)$ is called right $(H,\alpha,\bar{\sigma})$-Hom module.
\\

\textbf{Definition 4.4.} Let $(A,\alpha),(C,\beta)$ be monoidal Hom-algebras. Then $f: A\lr C$ is called weak Hom-algebra map if $f(1_{A})=1_{C}$, $f\circ\alpha=\beta\circ f$.
\\

\textbf{Definition 4.5.} Let $(H,\alpha)$ be a monoidal Hom-Hopf algebra and $(A,\beta)$ a monoidal Hom algebra. Let
 $\varphi^{+}: H\otimes A\lr A$ and $\varphi^{-}:A\otimes H\lr A$ are left and right actions, respectively.
  $(A,\beta)$ is called weak $(H,\alpha)$ Hom-bimodule if

(1) $\varphi^{+}(1_{H}\otimes a)=\beta(a)=\varphi^{-}(a\otimes 1_{H})$,

(2) $\varphi^{+}(\alpha\otimes\varphi^{-})=\varphi^{-}(\varphi^{+}\otimes\alpha)$.
\\

\textbf{Definition 4.6.} Let $(C_{\sharp_{\sigma}}^{\times}H, \beta\otimes\alpha)$ be monoidal Hom-bialgebra and $(A,\gamma)$ a monoidal Hom-bialgebra. Then $C\overset{p}\leftrightarrows_{j} A\overset{\pi}\rightleftarrows_{i} H$ is called a {\sl weak $m$-Hom admissible mappping system} if

(1) $p\circ j=id_{C}$, $\pi\circ i=id_{H}$,

(2) $\pi$ is a Hom-bialgebra map, $i$ is a weak Hom-algebra map and Hom-coalgebra map, $p$ is a Hom-coalgebra map, $j$ is a Hom-algebra map,

(3) $p$ is a weak $(H,\alpha)$ Hom-bimodule map [$(A,\gamma)$ is given a left action $ h\rightarrow a\triangleq i(\alpha^{-m}(h))a$ , a right action $a\leftarrow h\triangleq ai(\alpha^{-m}(h))$ by $(H,\alpha)$ and $(C,\beta)$ is given a right action $c\leftarrow h\triangleq \varepsilon(h)\beta(c)$ by $(H,\alpha)$]. And there is $\bar{\sigma}: H\otimes H\rightarrow A$ such that $(A,\gamma,\leftarrow), (A,\gamma,\rightarrow)$ are left, right $(H,\alpha,\bar{\sigma})$-Hom modules,

(4) $j(C)$ is a $(H,\alpha)$ Hom sub-bicomodule of $(A,\gamma)$ and $p|_{j(C)}$ is a Hom-bicomodule map. [$(A,\gamma)$ is given the $(H,\alpha)$ Hom-bicomodule defined by $\rho_{A}^{l}(a)=a_{(-1)}\otimes a_{(0)}\triangleq\alpha^{-m}\circ\pi(a_{1})\otimes a_{2}$ and $\rho_{A}^{r}(a)=a_{[0]}\otimes a_{[1]}\triangleq a_{1}\otimes \alpha^{-m}\circ\pi(a_{2})$; $(C,\beta)$ is given the $(H,\alpha)$ Hom-bicomodule defined by $\rho_{C}^{r}(c)=c_{[0]}\otimes c_{[1]}\triangleq\beta^{-1}(c)\otimes 1_{H}$],

(5) $(j\circ p)\ast(i\circ \pi)=id_{A}$.
\\

Now let $(C_{\sharp_{\sigma}}^{\times}H, \beta\otimes\alpha)$ be a monoidal Hom-bialgebra built in
 Theorem 3.3. Then there is some natural maps:
\begin{center}
  $\bar{j}: C_{\sharp_{\sigma}}^{\times}H\lr C, c\otimes h\mapsto \varepsilon(h)c$;\quad\quad
  $\bar{i}: C_{\sharp_{\sigma}}^{\times}H\lr H, c\otimes h\mapsto\varepsilon(c)h $;
\end{center}
\begin{center}
  $\bar{p}:C\lr C_{\sharp_{\sigma}}^{\times}H, c\mapsto c\otimes 1_{H}$;
  \quad\quad  $\bar{h}:H\lr C_{\sharp_{\sigma}}^{\times}H, h\mapsto 1_{C}\otimes h$.
\end{center}

A left, right $(H,\alpha)$ Hom-module action on $C_{\sharp_{\sigma}}^{\times}H$ defined by
\begin{center}
  $\varphi^{l}:l\otimes(a\otimes h )\mapsto(\alpha(l_{11})\rightarrow\beta^{-1}(a))\sigma(\alpha^{k+2-m}(l_{12}),\alpha^{k+1}(h_{1}))\otimes\alpha^{1-m}(l_{2})\alpha(h_{2})$
\end{center}
and
\begin{center}
  $\varphi^{r}:(a\otimes h)\otimes l\mapsto a\sigma(\alpha^{k+1}(h_{1}),\alpha^{k+1-m}(l_{1}))\otimes\alpha(h_{2})\alpha^{1-m}(l_{2})$.
\end{center}

A left,right $(H,\alpha)$ Hom-comodule structure on $C_{\sharp_{\sigma}}^{\times}H$ defined by
\begin{center}
  $\rho^{l}: a\otimes h\mapsto \alpha^{-1}(a_{(-1)})\alpha^{-1-m}(h_{1})\otimes a_{(0)}\otimes h_{2}$
\end{center}
and
\begin{center}
  $\rho^{r}: a\otimes h\mapsto \beta^{-1}(a)\otimes h_{1}\otimes\alpha^{-m}(h_{2})$.
\end{center}

With the above notion, we have

\textbf{Lemma 4.7.} If $(C_{\sharp_{\sigma}}^{\times}H, \beta\otimes\alpha)$ is a monoidal Hom-bialgebra, then $(C_{\sharp_{\sigma}}^{\times}H, \beta\otimes\alpha,\varphi^{l})$ is a left $(H,\alpha,\bar{\sigma})$ Hom module, where $\bar{\sigma}:h\otimes l\mapsto \sigma(\alpha^{k+1-m}(h),\alpha^{k+1-m}(l))\otimes 1_{H}$.

\textbf{Proof.} By Lemma 4.1, for any $a\otimes h\in C_{\sharp_{\sigma}}^{\times}H$, $l,g\in H$, we have
\begin{align*}
  &\varphi^{l}(\alpha\otimes\varphi^{l})[g\otimes l(a\otimes h)]\\
  =&\varphi^{l}[\alpha(g)\otimes(\alpha(l_{11})\rightarrow\beta^{-1}(a))\sigma(\alpha^{k+2-m}(l_{12}),\alpha^{k+1}(h_{1}))
  \otimes\alpha^{1-m}(l_{2})\alpha(h_{2})]\\
  =&[\alpha^{2}(g_{11})\rightarrow(l_{11}\rightarrow\beta^{-2}(a))\sigma(\alpha^{k+1-m}(l_{12}),\alpha^{k}(h_{1}))]\\
  &\sigma(\alpha^{k+3-m}(g_{12}),\alpha^{k+2-m}(l_{21})\alpha^{k+2}(h_{21}))
  \otimes \alpha^{2-m}(g_{2})[\alpha^{2-m}(l_{22})\alpha^{2}(h_{22})]\\
  =&[[\alpha^{2}(g_{111})\rightarrow(l_{11}\rightarrow\beta^{-2}(a))][\alpha^{2}(g_{112})\rightarrow
  \sigma(\alpha^{k+1-m}(l_{12}),\alpha^{k}(h_{1}))]]\\
  & \sigma(\alpha^{k+3-m}(g_{12}),\alpha^{k+2-m}(l_{21})
  \alpha^{k+2}(h_{21}))\otimes\alpha^{2-m}(g_{2})[\alpha^{2-m}(l_{22})\alpha^{2}(h_{22})]\\
  =&[\alpha(g_{11}l_{11})\rightarrow a][[\alpha^{2}(g_{121})\rightarrow\sigma(\alpha^{k+1-m}(l_{12}),\alpha^{k}(h_{1}))]\\
  &\sigma(\alpha^{k+3-m}(g_{122}),
  \alpha^{k+1-m}(l_{21})\alpha^{k+1}(h_{21}))]\otimes\alpha^{2-m}(g_{2})[\alpha^{2-m}(l_{22})\alpha^{2}(h_{22})]\\
 =&\{[\sigma(\alpha^{k+3-m}(g_{1111}),\alpha^{k+3-m}(l_{1111}))(\alpha(g_{1112}l_{1112})\rightarrow
  \beta^{-2}(a))]\sigma^{-1}(\alpha^{k+3-m}(g_{112}),\alpha^{k+3-m}(l_{112}))\}\\
  &[[\alpha^{2}(g_{121})\rightarrow
  \sigma(\alpha^{k+1-m}(l_{12}),\alpha^{k}(h_{1}))]\\
  & \sigma(\alpha^{k+3-m}(g_{122}),\alpha^{k+1-m}(l_{21})\alpha^{k+1}(h_{21}))]\otimes\alpha^{2-m}(g_{2})[\alpha^{2-m}(l_{22})\alpha^{2}(h_{22})]\\
  =&\{[\sigma(\alpha^{k+3-m}(g_{1111}),\alpha^{k+3-m}(l_{1111}))(\alpha(g_{1112}l_{1112})\rightarrow\beta^{-2}(a))]
  \sigma^{-1}(\alpha^{k+3-m}(g_{112}),\alpha^{k+3-m}(l_{112}))\}\\
  &\{[(\sigma(\alpha^{k+3-m}(g_{12111}),\alpha^{k+2-m}(l_{1211}))
  \sigma(\alpha^{k+2-m}(g_{12112})\alpha^{k+1-m}(l_{1212}),\alpha^{k}(h_{11})))\\
  &\sigma^{-1}(\alpha^{k+3-m}(g_{1212}),\alpha^{k+1-m}(l_{122})\alpha^{k}(h_{12}))]\sigma(\alpha^{k+3-m}(g_{122}),\alpha^{k+1-m}(l_{21})\alpha^{k+1}(h_{21})\}\\
  &\otimes\alpha^{2-m}(g_{2})[\alpha^{2-m}(l_{22})\alpha^{2}(h_{22})]\\
  =&[[\sigma(\alpha^{k+3-m}(g_{1111}),\alpha^{k+2-m}(l_{111}))(\alpha(g_{1112})l_{112}\rightarrow\beta^{-2}(a))]
  \sigma^{-1}(\alpha^{k+3-m}(g_{112}),\alpha^{k+2-m}(l_{12}))]\\
  &\{[\sigma(\alpha^{k+3-m}(g_{1211}),\alpha^{k+2-m}(l_{211}))\sigma(\alpha^{k+2-m}(g_{1212})\alpha^{k+1-m}(l_{212}),\alpha^{k}(h_{1}))]\\
  &[\sigma^{-1}(\alpha^{k+3-m}(g_{1221}),\alpha^{k+2-m}(l_{2211})\alpha^{k+1}(h_{211}))\\
  & \sigma(\alpha^{k+3-m}(g_{1222}),\alpha^{k+2-m}(l_{2212})
  \alpha^{k+1}(h_{212}))]\}\otimes\alpha^{2-m}(g_{2})[\alpha^{3-m}(l_{222})\alpha^{2}(h_{22})]\\
  =&[[\sigma(\alpha^{k+3-m}(g_{1111}),\alpha^{k+2-m}(l_{111}))(\alpha(g_{1112})l_{112}\rightarrow\beta^{-2}(a))]
  \sigma^{-1}(\alpha^{k+3-m}(g_{112}),\alpha^{k+2-m}(l_{12}))]\\
  &[\sigma(\alpha^{k+3-m}(g_{121}),\alpha^{k+2-m}(l_{211}))
  \sigma(\alpha^{k+2-m}(g_{122})\alpha^{k+1-m}(l_{212}),\alpha^{k}(h_{1}))]\\
  & \otimes\alpha^{2-m}(g_{2})[\alpha^{2-m}(l_{22})\alpha(h_{2})]\\
  =&[\sigma(\alpha^{k+3-m}(g_{111}),\alpha^{k+2-m}(l_{11}))(\alpha(g_{112})l_{12}\rightarrow\beta^{-1}(a))]
  [(\sigma^{-1}(\alpha^{k+3-m}(g_{1211}),\alpha^{k+3-m}(l_{2111}))\\
  &\sigma(\alpha^{k+3-m}(g_{1212}),\alpha^{k+3-m}(l_{2112})))\sigma(\alpha^{k+2-m}(g_{122})\alpha^{k+2-m}(l_{212}),
  \alpha^{k+1}(h_{1}))]\\
  & \otimes\alpha^{2-m}(g_{2})[\alpha^{2-m}(l_{22})\alpha(h_{2})]\\
  =&[\sigma(\alpha^{k+3-m}(g_{111}),\alpha^{k+2-m}(l_{11}))(\alpha(g_{112})l_{12}\rightarrow\beta^{-1}(a))]
  \sigma(\alpha^{k+2-m}(g_{12})\alpha^{k+2-m}(l_{21}),\alpha^{k+2}(h_{1}))\\
  &\otimes\alpha^{2-m}(g_{2})[\alpha^{2-m}(l_{22})\alpha(h_{2})]\\
  =&\sigma(\alpha^{k+3-m}(g_{11}),\alpha^{k+3-m}(l_{11}))[((g_{12}l_{12})\rightarrow\beta^{-1}(a))
  \sigma(\alpha^{k+1-m}(g_{21}l_{21}),\alpha^{k+1}(h_{1}))]\\
  &\otimes[\alpha^{2-m}(g_{22})\alpha^{2-m}(l_{22})]\alpha^{2}(h_{2})\\
  =&m_{C_{\sharp_{\sigma}}^{\times}H}[(\beta\otimes\alpha)\otimes\varphi^{l}](1\otimes m_{H}\otimes 1)
  (\bar{\rho}\otimes 1)[g\otimes l\otimes(a\otimes h)].
\end{align*}
This completes the proof.  $\hfill \blacksquare$
\\

\textbf{Lemma 4.8.} If $(C_{\sharp_{\sigma}}^{\times}H, \beta\otimes\alpha)$ is a monoidal Hom-bialgebra, then $(C_{\sharp_{\sigma}}^{\times}H, \beta\otimes\alpha,\varphi^{r})$ is a right $(H,\alpha,\bar{\sigma})$ Hom-module, where $\bar{\sigma}$ is as above.

\textbf{Proof.} For any $a\otimes h\in C_{\sharp_{\sigma}}^{\times}H$, $l,g\in H$, we have
\begin{align*}
  &\varphi^{r}(\varphi^{r}\otimes\alpha)[(a\otimes h)\otimes l\otimes g]=\varphi^{r}[a\sigma(\alpha^{k+1}(h_{1}),\alpha^{k+1-m}(l_{1}))\otimes\alpha(h_{2})\alpha^{1-m}(l_{2})\otimes\alpha(g)]\\
=&[a\sigma(\alpha^{k+1}(h_{1}),\alpha^{k+1-m}(l_{1}))]\sigma(\alpha^{k+2}(h_{21})\alpha^{k+2-m}(l_{21}),\alpha^{k+2-m}(g_{1}))\\
&\otimes[\alpha^{2}(h_{22})\alpha^{2-m}(l_{22})]\alpha^{2-m}(g_{2})\\
=&\beta(a)[\sigma(\alpha^{k+1}(h_{1}),\alpha^{k+1-m}(l_{1}))\sigma(\alpha^{k+1}(h_{21})\alpha^{k+1-m}(l_{21}),\alpha^{k+1-m}(g_{1}))]\\
&\otimes[\alpha^{2}(h_{22})\alpha^{2-m}(l_{22})]\alpha^{2-m}(g_{2})\\
=&\beta(a)[\sigma(\alpha^{k+2}(h_{11}),\alpha^{k+2-m}(l_{11}))\sigma(\alpha^{k+1}(h_{12})\alpha^{k+1-m}(l_{12}),\alpha^{k+1-m}(g_{1}))]\\
&\otimes[\alpha(h_{2})\alpha^{1-m}(l_{2})]\alpha^{2-m}(g_{2})\\
=&\beta(a)[[\alpha^{m+1}(h_{11})\rightarrow\sigma(\alpha^{k+1-m}(l_{11}),\alpha^{k+1-m}(g_{11}))]\sigma(\alpha^{k+2}(h_{12}),\alpha^{k+1-m}(l_{12}g_{12}))]
\\
&\otimes[\alpha(h_{2})\alpha^{1-m}(l_{2})]\alpha^{2-m}(g_{2})\\
=&m_{C_{\sharp_{\sigma}}^{\times}H}[(\beta(a)\otimes\alpha(h))\otimes(\sigma(\alpha^{k+2-m}(l_{1}),\alpha^{k+2-m}(g_{1}))\\
& \otimes\alpha^{1-m}(l_{2}g_{2}))]\\
=&m_{C_{\sharp_{\sigma}}^{\times}H}[(\beta\otimes\alpha)\otimes\varphi^{r}][(a\otimes h)\otimes[(\sigma(\alpha^{k+1-m}(l_{1}),\alpha^{k+1-m}(g_{1}))\otimes1)\otimes l_{2}g_{2}]]\\
=&m_{C_{\sharp_{\sigma}}^{\times}H}[(\beta\otimes\alpha)\otimes\varphi^{r}][(a\otimes h)\otimes\bar{\sigma}(l_{1},g_{1})\otimes l_{2}g_{2}]\\
=&m_{C_{\sharp_{\sigma}}^{\times}H}[(\beta\otimes\alpha)\otimes\varphi^{r}][1\otimes1\otimes m_{H}](1\otimes \bar{\rho})[(a\otimes h)\otimes l\otimes g]
\end{align*}
This completes the proof. $\hfill \blacksquare$
\\

\textbf{Lemma 4.9.} If $(C_{\sharp_{\sigma}}^{\times}H, \beta\otimes\alpha)$ is a monoidal Hom-bialgebra, then $(C_{\sharp_{\sigma}}^{\times}H, \beta\otimes\alpha,\varphi^{l},\varphi^{r})$ is a weak $(H,\alpha)$ Hom-bimodule.

\textbf{Proof.} For any $a\otimes h\in C_{\sharp_{\sigma}}^{\times}H$, $l,g\in H$, we have
\begin{align*}
  &\varphi^{l}(\alpha\otimes\varphi^{r})[l\otimes (a\otimes h)\otimes g]=\varphi^{l}[\alpha(l)\otimes a\sigma(\alpha^{k+1}(h_{1}),\alpha^{k+1-m}(g_{1}))\otimes\alpha(h_{2})\alpha^{1-m}(g_{2})]\\
  =&[\alpha^{2}(l_{11})\rightarrow\beta^{-1}(a)\sigma(\alpha^{k}(h_{1}),\alpha^{k-m}(g_{1}))]\sigma(\alpha^{k+3-m}(l_{12}),\alpha^{k+2}(h_{21})\alpha^{k+2-m}(g_{21}))\\&\otimes\alpha^{2-m}(l_{2})[\alpha^{2}(h_{22})\alpha^{2-m}(g_{22})]\\
  =&[(\alpha^{2}(l_{111})\rightarrow\beta^{-1}(a))(\alpha^{2}(l_{112})\rightarrow\sigma(\alpha^{k}(h_{1}),\alpha^{k-m}(g_{1})))]\sigma(\alpha^{k+3-m}(l_{12}),\alpha^{k+2}(h_{21})\alpha^{k+2-m}(g_{21}))\\&\otimes\alpha^{2-m}(l_{2})[\alpha^{2}(h_{22})\alpha^{2-m}(g_{22})]\\
  =&(\alpha^{2}(l_{11})\rightarrow a)[(\alpha^{2}(l_{121})\rightarrow\sigma(\alpha^{k+1}(h_{11}),\alpha^{k+1-m}(g_{11})))\sigma(\alpha^{k+3-m}(l),\alpha^{k+1}(h_{12})\alpha^{k+1-m}(g_{12}))]\\&\otimes\alpha^{2-m}(l_{2})[\alpha(h_{2})\alpha^{1-m}(g_{2})]\\
  =&(\alpha^{2}(l_{11})\rightarrow a)[\sigma(\alpha^{k+3-m}(l_{121}),\alpha^{k+2}(h_{11}))\sigma(\alpha^{k+2-m}(l_{122})\alpha^{k+1}(h_{12}),\alpha^{k+1-m}(g_{1}))]\\&\otimes\alpha^{2-m}(l_{2})[\alpha(h_{2})\alpha^{1-m}(g_{2})]\\
  =&[(\alpha(l_{11})\rightarrow \beta^{-1}(a))\sigma(\alpha^{k+2-m}(l_{12}),\alpha^{k+1}(h_{1}))]\sigma(\alpha^{k+2-m}(l_{21})\alpha^{k+2}(h_{21}),\alpha^{k+2-m}(g_{1}))\\&\otimes[\alpha^{2-m}(l_{22})\alpha^{2}(h_{22})]\alpha^{2-m}(g_{2})\\
  =&\varphi^{r}[(\alpha(l_{11})\rightarrow\beta^{-1}(a))\sigma(\alpha^{k+2-m}(l_{12}),\alpha^{k+1}(h_{1}))\otimes\alpha^{1-m}(l_{2})\alpha(h_{2})\otimes\alpha(g)]\\
  =&\varphi^{r}(\varphi^{l}\otimes\alpha)[l\otimes (a\otimes h)\otimes g]
  \end{align*}
  This completes the proof. $\hfill \blacksquare$
\\

\textbf{Proposition 4.10.} If $(C_{\sharp_{\sigma}}^{\times}H, \beta\otimes\alpha)$ is a monoidal Hom-bialgebra built in Theorem 3.3, then $C\overset{\bar{j}}\leftrightarrows_{\bar{p}} C_{\sharp_{\sigma}}^{\times}H\overset{\bar{i}}\rightleftarrows_{\bar{\pi}} H$ is a weak $m$-Hom
admissible mappping system.
\\

\textbf{Theorem 4.11.} Let $C\overset{p}\leftrightarrows_{j} A\overset{\pi}\rightleftarrows_{i} H$ be a weak $m$-Hom admissible mappping system.
 Then  $C_{\sharp_{\sigma}}^{\times}H\cong A$ as monoidal Hom-bialgebra.

\textbf{Proof.} Define the following maps:
\begin{eqnarray*}
   &f: C_{\sharp_{\sigma}}^{\times}H\rightarrow A, c\otimes h\mapsto\gamma^{-1}(j(c)i(h))  \\
    &g:A\rightarrow C_{\sharp_{\sigma}}^{\times}H,  a\mapsto(\beta\otimes\alpha)(p(a_{1})\otimes \pi(a_{2})).
     \end{eqnarray*}
We need to prove the following aspects:
\begin{itemize}
  \item $f\circ g=id_{A}$, $g\circ f=id_{C_{\sharp_{\sigma}}^{\times}H}$;
  \item $f$ is a Hom algebra homomorphism;
  \item $g$ is a Hom coalgebra homomorphism.
\end{itemize}

Firstly, it is easy to verify $f\circ(\beta\otimes\alpha)=\gamma\circ f$, $g\circ\gamma=(\beta\otimes\alpha)\circ g$, then for any $a\in A$,
\begin{align*}
  f(g(a))=&f[(\beta\otimes\alpha)(p(a_{1})\otimes \pi(a_{2}))]\\
  =&\gamma\circ f[(p(a_{1})\otimes \pi(a_{2}))]\\
  =&j(p(a_{1}))i(\pi(a_{2}))=(j\circ p)\ast(i\circ\pi)(a)=a.
\end{align*}
And
\begin{align*}
  g\circ f(c\otimes h)&=g[\gamma^{-1}(j(c)i(h))]\\
  =&(\beta^{-1}\otimes\alpha^{-1})g[(j(c)i(h))]\\
  =&p((j(c)_{1}i(h)_{1}))\otimes\pi((j(c)_{2}i(h)_{2}))\\
  =&p((j(c)_{1}i(h_{1})))\otimes\pi((j(c)_{2}i(h_{2})))\\
  =&\varepsilon(h_{1})\beta[p(j(c)_{1})]\otimes\pi[j(c)_{2}]h_{2}\\
  =&\beta[\beta^{-1}(c)]\otimes 1_{H}\alpha^{-1}(h_{2})=c\otimes h.
\end{align*}

Next, we prove $f$ is a Hom algebra homomorphism, we only need to prove $f[(a\otimes h)(b\otimes g)]=f(a\otimes h)f(b\otimes g)$.
First we compute for any $h\in H$, $b\in C$
\begin{align*}
  i(h)j(b)=&j\circ p[i(h_{1})j(b)_{1}]i\circ\pi[i(h_{2})j(b)_{2}]\\
  =&j\circ p[\alpha^{m}(h_{1})\rightarrow j(b)_{1}]i[h_{2}\pi(j(b)_{2})]\\
  =&j[\alpha^{m}(h_{1})\rightarrow p(j(b)_{1})]i[h_{2}\pi(j(b)_{2})]\\
  =&j[\alpha^{m}(h_{1})\rightarrow \beta^{-1}(b)]i[\alpha(h_{2})]
\end{align*}
Therefore,
\begin{align*}
  & f[(a\otimes h)(b\otimes g)]\\
  =&\gamma^{-1}\{[j(a)[j(\alpha^{m}(h_{11})\rightarrow\beta^{-2}(b))j(\sigma(\alpha^{k+1}(h_{12}),\alpha^{k}(g_{1})))]]i(\alpha(h_{2}g_{2}))\}\\
  =&\gamma^{-1}\{[j(a)j(\alpha^{m}(h_{1})\rightarrow\beta^{-1}(b))][j(\sigma(\alpha^{k+2}(h_{21}),\alpha^{k+1}(g_{1})))i(\alpha(h_{22})g_{2})]\}\\
  =&\gamma^{-1}\{[j(a)j(\alpha^{m}(h_{1})\rightarrow\beta^{-1}(b))]\beta[j(\sigma(\alpha^{k+1}(h_{21}),\alpha^{k}(g_{1})))\leftarrow\alpha^{m}(h_{22})\alpha^{m-1}(g_{2})]\}\\
  =&\gamma^{-1}\{[j(a)j(\alpha^{m}(h_{1})\rightarrow\beta^{-1}(b))]\beta[\bar{\sigma}(\alpha^{m}(h_{21}),\alpha^{m-1}(g_{1}))\leftarrow\alpha^{m}(h_{22})\alpha^{m-1}(g_{2})]\}\\
  =&\gamma^{-1}\{[j(a)j(\alpha^{m}(h_{1})\rightarrow\beta^{-1}(b))][(1_{A}\leftarrow\alpha^{m}(h_{2}))\leftarrow\alpha^{m}(g)]\}\\
  =&\gamma^{-1}\{[j(a)j(\alpha^{m}(h_{1})\rightarrow\beta^{-1}(b))][i(\alpha(h_{2}))i(g)]\}\\
  =&\gamma^{-1}\{j(\alpha(a))[[j(\alpha^{m-1}(h_{1})\rightarrow\beta^{-2}(b))i(h_{2})]i(g)]\}\\
  =&\gamma^{-1}\{j(\alpha(a))[[i(\alpha^{-1}(h))j(\beta^{-1}(b))]i(g)]\}\\
  =&\gamma^{-1}\{[j(a)i(h)][j(b)i(g)]\}\\
  =&f(a\otimes h)f(b\otimes g).
\end{align*}

Finally, we prove $g$ is a Hom coalgebra homomorphism, we only need to prove $\Delta(g(a))=(g\otimes g)\Delta(a)$. First by Def 4.6.(3)and(4), it is not hard to verify for any $a\in A$,
\begin{equation}
  \alpha^{m+1}(p(a_{1})_{(-1)})\pi(a_{2})\otimes\beta(p(a_{1})_{(0)})=\alpha\circ\pi(a_{1})\otimes p(a_{2}).
\end{equation}
Therefore,
\begin{align*}
  \Delta(g(a))=&p(\alpha(a_{11}))\otimes\alpha^{m+1}[p(a_{12})_{(-1)}]\pi(a_{21})\otimes\beta^{2}[p(a_{12})(0)]\otimes\alpha[\pi(a_{22})]\\
  =&p(a_{1})\otimes\alpha^{m+2}[p(a_{211})_{(-1)}]\alpha[\pi(a_{212})]\otimes\beta^{3}[p(a_{211})(0)]\otimes\alpha[\pi(a_{22})]\\
  \overset{(4.1)}=&p(a_{1})\otimes\alpha^{2}[\pi(a_{211})]\otimes\beta^{2}[p(a_{212})]\otimes\alpha[\pi(a_{22})]\\
  =&(g\otimes g)\Delta(a)
\end{align*}

In fact, it is easy to prove $f$ is a Hom coalgebra homomorphism and $g$ is a Hom coalgebra homomorphism by the relation
of $f$ and $g$. Thus, $f$ and $g$ are Hom bialgebra homomorphisms and
 $C_{\sharp_{\sigma}}^{\times}H\cong A$.

 This completes the proof. $\hfill \blacksquare$

{\bf Remark 4.12.} In the case of Hopf algebras,  it follows  from Theorem 4.11 that \cite[Theorem 3.9]{WJZ}
  when we take  $m=0$, $k=-1$ and $\sigma$ is trivial.

\section*{Acknowledgements}

The second author thanks the financial support of the National Natural Science Foundation of
 China (Grant No. 11871144 and No. 12271089).

\end{document}